\documentclass{amsart}
\usepackage{bbm,amsmath,amsfonts,amssymb}

\newtheorem{theorem}{Theorem}[section]
\newtheorem{lemma}[theorem]{Lemma}

\newtheorem{proposition}[theorem]{Proposition}
\newtheorem{definition}[theorem]{Definition}

\newtheorem{remark}[theorem]{Remark}

\newcounter{figures}[section]

\def\bZ{{\mathbb Z}}

\def\bN{{\mathbb N}}

\def\bR{{\mathbb R}}

\def\cB{\mathcal{B}}

\def\cF{\mathcal{F}}

\def\cJ{\mathcal{J}}

\def\cM{\mathcal{M}}
\def\cN{\mathcal{N}}

\def\cP{\mathcal{P}}

\def\cS{\mathcal{S}}

\def\LL{L}
\def\supp{{\rm{supp}\, }}
\def\eps{{\varepsilon}}

\def\diam{{\rm diam \,}}
\def\ph{\varphi}

\def\R{\mathbb{R}}

\def\cs{{c_\sharp}}
\def\ct{\tilde{c}}

\def\Cstar{{C^\star}}
\def\cstar{{c^\star}}

\def\cc{\tau}
\def\kr{{s_r}}
\def\dist{{\rm{dist}}}
\def\kkl{{k_0}}
\def\kl{{\kappa_0}}
\def\kll{{\kappa_1}}
\def\ddelta{t}
\def\NN{K}

\begin{document}

\title
{Hardy spaces associated with non-negative self-adjoint operators}



\author{S. Dekel}
\address{Hamanofim St. 9, Herzelia, Israel}
\email{Shai.Dekel@ge.com}

\author{G. Kerkyacharian}
\address{
Laboratoire de Probabilit\'{e}s et Mod\`{e}les Al\'{e}atoires, CNRS-UMR 7599,
Universit\'{e} Paris VI et Universit\'{e} Paris VII, rue de Clisson, F-75013 Paris}
\email{kerk@math.univ-paris-diderot.fr}

\author{G. Kyriazis}
\address{Department of Mathematics and Statistics,
University of Cyprus, 1678 Nicosia, Cyprus}
\email{kyriazis@ucy.ac.cy}

\author{P. Petrushev}
\address{Department of Mathematics\\University of South Carolina\\
Columbia, SC 29208}
\email{pencho@math.sc.edu}




\begin{abstract}
Maximal and atomic Hardy spaces $H^p$ and $H^p_A$, $0<p\le 1$, are considered
in the setting of a doubling metric measure space in the presence of
a non-negative self-adjoint operator whose heat kernel has Gaussian localization
and the Markov property.
It is shown that $H^p=H^p_A$ with equivalent norms.
\end{abstract}


\maketitle

\section{Introduction}\label{sec:introduction}
\setcounter{equation}{0}

The purpose of this article is to establish the equivalence of
the maximal and atomic Hardy spaces $H^p$ and $H^p_A$, $0<p\le 1$,
in the general setting of a metric measure space with the doubling property
and in the presence of a non-negative self-adjoint operator whose heat kernel has Gaussian localization
and the Markov property.
%
We next describe our setting in detail (see \cite{CKP, KP}):


I. We assume that $(M,\rho,\mu)$ is a metric measure space satisfying the conditions:
$(M, \rho)$ is a locally compact  metric space with distance $\rho(\cdot, \cdot)$
and $\mu$ is a~positive Radon measure
such that the following {\em volume doubling condition} is valid
\begin{equation}\label{doubling-0}
0 < \mu(B(x,2r) ) \le c_0\mu(B(x,r))<\infty
\quad\hbox{for all $x \in M$ and $r>0$,}
\end{equation}
where $B(x,r)$ is the open ball centered at $x$ of radius $r$ and $c_0>1$ is a constant.
From above it follows that
\begin{equation}\label{doubling}
\mu(B(x,\lambda r) ) \le c_0\lambda^d \mu(B(x,r))
\quad\hbox{for $x \in M$, $r>0$, and $\lambda >1$,}
\end{equation}
were $d=\log_2 c_0 >0$ is a constant playing the role of a dimension.


II. The main assumption is that the local geometry of the space $(M,\rho,\mu)$ is related to
an essentially self-adjoint non-negative operator $L$ on $L^2(M, d\mu)$,
mapping real-valued to real-valued functions,
such that the associated semigroup $P_t=e^{-tL}$ consists of integral operators with
(heat) kernel $p_t(x,y)$ obeying the conditions:

\smallskip

\noindent
(a) Gaussian upper bound:
\begin{equation}\label{Gauss-local}
|p_t(x,y)|
\le \frac{ \Cstar\exp\{-\frac{\cstar\rho^2(x,y)}t\}}{\sqrt{\mu(B(x,\sqrt t))\mu(B(y,\sqrt t))}}
\quad\hbox{for} \;\;x,y\in M,\,t>0.
\end{equation}

\noindent
(b) H\"{o}lder continuity: There exists a constant $\alpha>0$ such that
\begin{equation}\label{lip}
\big|  p_t(x,y) - p_t(x,y')  \big|
\le \Cstar\Big(\frac{\rho(y,y')}{\sqrt t}\Big)^\alpha
\frac{\exp\{-\frac{\cstar\rho^2(x,y)}t \}}{\sqrt{\mu(B(x,\sqrt t))\mu(B(y, \sqrt t))}}
\end{equation}
for $x, y, y'\in M$ and $t>0$, whenever $\rho(y,y')\le \sqrt{t}$.

\smallskip

\noindent
(c) Markov property:
\begin{equation}\label{hol3}
\int_M p_t(x,y) d\mu(y)= 1
\quad\hbox{for $x\in M$ and $t >0$.}
\end{equation}
Above $\Cstar, \cstar>0$ are structural constants.

We also stipulate the following conditions on the geometry of $M$:

\noindent
(d) Reverse doubling condition: There exists a constant
$c_1>1$ such that
\begin{equation}\label{reverse-doubling}
\mu(B(x,2r) ) \ge c_1 \mu(B(x,r))
\quad\hbox{for $x \in M$ and $0< r \le \frac{\diam M}{3}$.}
\end{equation}

\smallskip

\noindent
(e) Non-collapsing condition in the case when $\mu (M) =\infty$:
There exists a~constant $c_2>0$ such that
\begin{equation}\label{non-collapsing}
\inf_{x\in M}\mu(B(x,1) )\ge c_2.
\end{equation}

A natural effective realization of the above setting appears in the general
framework of strictly local regular Dirichlet spaces with a complete intrinsic metric,
where it suffices to only verify the local Poincar\'e inequality and the global doubling condition
on the measure and then the above general setting applies in full.
In~particular, this setting covers the cases of
Lie groups or homogeneous spaces with polynomial volume growth,
complete Riemannian manifolds with Ricci curvature bounded from below and satisfying
the volume doubling condition.
Naturally, it contains the classical setting on $\R^n$
as well as the cases of the sphere, interval, ball, and simplex with weights.
For more details, see \cite{CKP}.

The maximal Hardy space $H^p$, $0<p\le 1$, will be defined (Definition~\ref{def:Hp})
as a set of distributions in the above described setting via the quasi-norm
$$
\|f\|_{H^p}:= \big\|\sup_{t>0}|e^{-t^2L}f(\cdot)|\big\|_{L^p}
$$
and its equivalence with a quasi-norm defined by the respective grand maximal operator will be established.
Following to some extent \cite{HLMMY, Duong-Li} the atomic Hardy space $H^p_A$, $0<p\le 1$,
will be defined (Definition~\ref{def-Hardy-Atomic-hom}) via atoms $a(x)$ with the properties:
There exists a function $b\in D(L^n)$ and a ball $B$ of radius $r=r_B >0$
such that
\begin{enumerate}
\item[(i)] $a=L^nb,$
\item[(ii)] $\supp L^k b\subset B$, $k=0, 1, \dots, n$, and
\item[(iii)] $\|L^k b\|_\infty \le r^{2(n-k)}|B|^{-1/p}$, $k=0, 1, \dots, n$,
\end{enumerate}
with $n:=\lfloor d/2p\rfloor +1$.
Naturally, an additional kind of atoms supported on $M$ will be introduced in the compact case.

Our main theorem (Theorem~\ref{thm:hardy}) asserts that in the setting described above
the maximal and atomic Hardy spaces $H^p$ and $H^p_A$ are the same
for $0<p\le 1$ with equivalent quasi-norms.
To prove the nontrivial embedding $H^p\subset H^p_A$ we devise a completely new approach
to atomic decomposition of Hardy spaces.
Our method relies on an idea different from the one of the classical proof, in particular,
it does not use the Calder\'{o}n-Zygmund decomposition.
This method is new when applied in the classical setting on $\R^n$ as well.
In light of the development of Besov and Triebel-Lizorkin spaces from \cite{KP},
this result shows that although general our setting allows to develop
the Littlewood-Paley theory and function spaces in almost complete analogy with the classical
case on $\R^n$.

This paper is organized as follows.
In \S\ref{sec:background} we assemble the necessary background material from \cite{CKP, KP}.
In \S\ref{sec:Hp-spaces} we introduce the maximal Hardy spaces and establish their characterization via
several maximal operators.
In \S\ref{sec:atomic-Hp-spaces} we introduce the atomic Hardy spaces and
in \S\ref{sec:equivalence-Hp-spaces} we prove our main result: the equivalence of maximal and atomic Hardy spaces.
In \S\ref{sec:littlewood-paley} we briefly discuss the characterization of Hardy spaces via square functions.
Section \ref{sec:appendix-1} is an appendix where we place the proofs of some auxiliary
assertions from previous sections.

\smallskip

\noindent
{\bf Notation.} For any set $E\subset M$ and $x\in M$ we denote
$\dist (x, E):= \inf_{y\in E} \rho(x, y)$, $E^c:= M\setminus E$, and $|E|:= \mu(E)$.
We will use the notation
$cB(x, \delta):= B(x, c\delta)$.
The class of Schwartz functions on $\R$ will be denoted by $\cS(\R)$.
As usual $C^\infty_0(\R)$ will stand for the class of all compactly supported $C^\infty$ functions on $\R$
and $C(E)$ will be the set of all continuous functions on $E$.
Positive constants will be denoted by $c$, $c_1$, $c'$, $\dots$ and they may vary
at every occurrence. Most of them will depend on the basic structural constants $c_0, \Cstar, \cstar$
from (\ref{doubling-0})-(\ref{lip}). This dependence usually will not be indicated explicitly.
The notation $a\sim b$ will mean $c_1\le a/b\le c_2$.

\section{Background}\label{sec:background}
\setcounter{equation}{0}

Our development of Hardy spaces will rely on some basic facts and results from \cite{CKP, KP},
which we review next.
We begin with the observation that as $L$ is a non-negative self-adjoint operator
that maps real-valued to real-valued functions, then for any real-valued, measurable and
bounded function $f$ on $\R_+$ the operator $f(L)$, defined by
$f(L):=\int_0^\infty f(\lambda)dE_\lambda$
with $E_\lambda$, $\lambda \ge 0$, being the spectral resolution associated with $L$,
is bounded on $L^2$, self-adjoint, and maps real-valued functions to real-valued functions.
Furthermore, if $f(L)$ is an integral operator, then its kernel $f(L)(x, y)$ is real-valued and
$f(L)(y, x)=f(L)(x, y)$,
in particular, $p_t(x, y)\in \R$ and $p_t(y,x) = p_t(x, y)$.

\subsection{Functional calculus}\label{subsec:func-calculus}

The {\em finite speed propagation property} plays a crucial role in this theory:
\begin{equation}\label{finite-speed}
\big\langle \cos(t\sqrt{L})f_1, f_2 \big\rangle=0,
\quad 0< \ct t<r,
\quad \ct:=\frac{1}{2\sqrt{\cstar}},
\end{equation}
for all open sets $U_j \subset M$, $f_j\in \LL^2(M)$, $\supp f_j\subset U_j$,
$j=1, 2$, where $r:=\rho(U_1, U_2)$.

This property leads to the following localization result for
the kernels of operators of the form $f(t\sqrt{L})$ whenever $\hat f$ is band limited.
Here $\hat f(\xi):=\int_\R f(x)e^{-ix\xi}dx$.


\begin{proposition}\label{prop:finite-sp}
Let $f$ be even, $\supp \hat f \subset [-A, A]$ for some $A>0$,
and $\hat f\in W^m_1$ for some $m>d$, i.e. $\|\hat f^{(m)}\|_{L^1} <\infty$.
Then for any $\ddelta>0$ and $x, y\in M$
\begin{equation}\label{finite-speed-2}
f(\ddelta\sqrt{L})(x, y) = 0
\quad\hbox{if}\quad
\rho(x, y) > \ct \ddelta A.
\end{equation}
\end{proposition}

We will need the following result which follows from \cite[Theorem 3.4]{KP} and (\ref{D2}).


\begin{theorem}\label{thm:S-local-kernels}
Suppose $f\in C^m(\bR_+)$, $m\ge d+1$,
$$
\hbox{
$|f^{(\nu)}(\lambda)|\le A_m(1+\lambda)^{-r}$ for $\lambda\ge 0$ and $0\le \nu\le m$, where $r > m+d$,
}
$$
and $f^{(2\nu+1)}(0)=0$ for $\nu\ge 0$ such that $2\nu+1 \le m$.
Then $f(\ddelta \sqrt L)$ is an integral operator with kernel $f(\ddelta \sqrt L)(x, y)$
satisfying
$$ 
\big|f(\ddelta \sqrt L)(x, y)\big|
\le cA_m |B(x, \ddelta)|^{-1}\big(1+\ddelta^{-1}\rho(x, y)\big)^{-m+d/2}
\quad\hbox{and}
$$ 
$$ 
\big|f(\ddelta \sqrt L)(x, y)-f(\ddelta \sqrt L)(x, y')\big|
\le cA_m |B(x, \ddelta)|^{-1}\Big(\frac{\rho(y, y')}{\ddelta}\Big)^\alpha
\big(1+\ddelta^{-1}\rho(x, y)\big)^{-m+d/2},
$$ 
whenever $\rho(y, y')\le t$.
Here $\alpha>0$ is from $(\ref{lip})$
and 
$c>0$ is a~constant depending only on $r, m$ and
the structural constants $c_0, \Cstar, \cstar$, $\alpha$.

Moreover, $\int_M f(\ddelta \sqrt L)(x, y)d\mu(y)=f(0)$.
\end{theorem}


\begin{remark}\label{rem:background}
Observe that Theorem~\ref{thm:S-local-kernels} is established in \cite[Theorem 3.4]{KP}
in the case when $0<t\le 1$. However the same proof applies also to the case $0<t <\infty$.
\end{remark}

\subsection{On the geometry of the underlying space}\label{subsec:geometry}

As is shown in \cite[Proposition 2.2]{CKP} if $M$ is connected the reverse doubling condition
(\ref{reverse-doubling}) follows from the doubling condition (\ref{doubling-0}) and
hence it is not very restrictive.
Note that (\ref{reverse-doubling}) implies
\begin{equation}\label{reverse-doubling-2}
|B(x,\lambda r)| \ge c\lambda^\eps |B(x, r)|,
\quad \lambda > 1, \; r>0,\; \hbox{$0< \lambda r\le \frac{\diam M}{3}$,}
\end{equation}
where $c, \eps >0$ are constants.
This coupled with (\ref{non-collapsing}) leads to
\begin{equation}\label{ball-low-est}
\inf_{x\in M}|B(x, r)| \ge cr^\eps \quad\hbox{for}\quad r>1.
\end{equation}
In \cite[Proposition 2.1]{CKP} it is shown that
$\mu(M) < \infty$ if and only if $\diam M<\infty$.
Then denoting $D:=\diam M$ we obtain using (\ref{doubling})
\begin{equation}\label{ball-low-est-2}
\inf_{x\in M} |B(x, r)| \ge r^d c_0^{-1}D^{-d}\mu(M), \quad 0<r\le D,
\end{equation}
which is a substitute for (\ref{ball-low-est}) in the case when $\mu(M)<\infty$.

\smallskip

To compare the volumes of balls with different centers $x, y\in M$ and the same radius $r$
we will use the inequality
\begin{equation}\label{D2}
|B(x, r)| \le c_0\Big(1+ \frac{\rho(x,y)}{r}\Big)^d  |B(y, r)|,
\quad x, y\in M, \; r>0.
\end{equation}
As $B(x,r) \subset B(y, \rho(y,x) +r)$ the above inequality is immediate from (\ref{doubling}).

\smallskip

The following simple inequalities will also be needed \cite[Lemma 2.1]{KP}:
For $\sigma >d$ and $\ddelta>0$
\begin{equation}\label{tech-1}
\int_M \big(1+\ddelta^{-1}\rho(x, y)\big)^{-\sigma} d\mu(y) \le c|B(x, \ddelta)|,
\quad x\in M,
\end{equation}
\begin{equation}\label{tech-2}
\int_M \big(1+\ddelta^{-1}\rho(x, u)\big)^{-\sigma} \big(1+\ddelta^{-1}\rho(u, y)\big)^{-\sigma}d\mu(u)
\le c|B(x, \ddelta)|\big(1+\ddelta^{-1}\rho(x, y)\big)^{-\sigma+d}.
\end{equation}

\subsection{Distributions}\label{sybsec:distributions}

The Hardy spaces $H^p$, $0<p\le 1$, associated with $L$ will be spaces of distributions.
In the setting of this article the class of test functions $\cS= \cS(L)$
is defined (see \cite{KP})  as the set of all complex-valued functions
$\phi\in  \cap_{m\ge 1} D(L^m)$ such that
\begin{equation}\label{norm-S}
\cP_{m}(\phi)
:= \sup_{x\in M} (1+\rho(x, x_0))^m \max_{0\le \nu\le m}|L^\nu\phi(x)| < \infty,
\quad \forall m \ge 0.
\end{equation}
Here $x_0\in M$ is selected arbitrarily and fixed once and for all.
Observe that if $\phi\in\cS$, then $\overline{\phi} \in \cS$,
which is a consequence of the fact that $L \overline{\phi}= \overline{L\phi}$,
for $L$ maps real-valued to real-valued functions.

As usual the space $\cS'$ of distributions on $M$ is defined as the set of all
continuous linear functionals on $\cS$
and the action of $f\in \cS'$ on $\overline{\phi}\in \cS$ will be denoted by
$\langle f, \phi\rangle:= f(\overline{\phi})$,
which is consistent with the inner product on $L^2(M)$.
Clearly, for any $f\in \cS'$ there exist constants $m\in \bZ_+$ and $c>0$ such that
\begin{equation}\label{distribution-1}
|\langle f, \phi\rangle| \le c\cP_m(\phi), \quad \forall \phi\in\cS.
\end{equation}

It will be useful to clarify the action of operators of the form $\varphi(\sqrt L)$ on  $\cS'$.
Note first that if the function $\varphi\in \cS(\R)$ is real-valued and even, then
from Theorem~\ref{thm:S-local-kernels} it follows that
$\varphi(\sqrt{L})(x, \cdot)\in \cS$
and $\varphi(\sqrt{L})(\cdot, y)\in \cS$.
Furthermore, it is easy to see that  $\varphi(\sqrt L)$ maps continuously $\cS$  into $\cS$.


\begin{definition}\label{def:phi-distr}
We define $\varphi(\sqrt L)f$ for any $f\in \cS'$ by
\begin{equation}\label{phi-distr}
\langle \varphi(\sqrt L)f, \phi \rangle := \langle f, \varphi(\sqrt L)\phi \rangle
\quad\hbox{for}\quad \phi \in \cS.
\end{equation}
\end{definition}

From above it follows that, $\varphi(\sqrt L)$ maps continuously $\cS'$  into $\cS'$.
Furthermore, if $\varphi, \psi\in \cS(\R)$ are real-valued and even, then
\begin{equation}\label{comute}
\varphi(\sqrt{L})\psi(\sqrt{L})f =  \psi(\sqrt{L})\varphi(\sqrt{L}) f, \quad \forall f\in\cS'.
\end{equation}


\begin{proposition}\label{prop:dist}
Suppose $\varphi\in\cS$ is real-valued and even and let $f\in \cS'$.
Then
\begin{equation}\label{phif}
\varphi(\sqrt L)f(x) = \langle  f , \varphi(\sqrt L)(x,\cdot)\rangle, \quad x\in M.
\end{equation}
Moreover, $\varphi(\sqrt L)f$ is a continuous and slowly growing function,
namely, there exist constants $m\in\bZ_+$ and $c>0$, depending on $f$, such that
\begin{equation}\label{slowly}
|\varphi(\sqrt L)f(x)| \le c(1+\rho(x, x_0))^m, \quad x\in M, \quad\hbox{and}
\end{equation}
\begin{equation}\label{lipsch}
|\varphi(\sqrt L)f(x)-\varphi(\sqrt L)f(x')| \le c(1+\rho(x, x_0))^m \rho(x, x')^\alpha,
\quad\hbox{if}\;\; \rho(x, x')\le 1.
\end{equation}
Here $\alpha>0$ is the constant from $(\ref{lip})$.
\end{proposition}

To streamline our exposition we place the proof of this assertion in the appendix.

We now give the main convergence result for distributions.


\begin{proposition}\label{prop:convergence}
Suppose $\varphi \in \cS(\R)$, $\varphi$ is real-valued and even, and $\varphi(0)=1$.
Then for every $\phi\in \cS$
\begin{equation}\label{decomp-dist-1}
\phi=\lim_{\ddelta\to 0} \varphi (\ddelta\sqrt L) \phi
\quad \mbox{ in }\; \cS,
\end{equation}
and for every $f\in \cS'$
\begin{equation}\label{decomp-dist-2}
f=\lim_{\ddelta\to 0} \varphi (\ddelta\sqrt L) f
\quad \mbox{ in }\; \cS'.
\end{equation}
Furthermore,
for any $f\in L^p(M)$, $1\le p < \infty$, $(\ref{decomp-dist-2})$ is valid with the convergence in $L^p$.
In addition, for any $f\in L^p(M)$, $1\le p\le \infty$, one has
$f(x)=\lim_{\ddelta\to 0} \varphi (\ddelta\sqrt L) f(x)$ for almost all $x\in M$.

\end{proposition}

This proposition is established in \cite[Proposition 5.5]{KP}
in the case when $\varphi$ is compactly supported and $\varphi^{(\nu)}(0)=0$ for $\nu\ge 1$.
We give its proof in the appendix.

For more information on distributions in the setting of this paper, see \cite{KP}.

\section{Hardy spaces via maximal operators}\label{sec:Hp-spaces}
\setcounter{equation}{0}

In this section we introduce several maximal operators and define the Hardy spaces $H^p$, $0<p\le 1$,
in the setting described in the introduction.
As in the classical case on $\R^n$ the grand maximal operator will play an important r\^{o}le.

\subsection{Maximal operators and definition of \boldmath $H^p$}\label{sec:max-operators}


\begin{definition}\label{def:admissible}
A function $\varphi \in\cS(\R)$ is called {\em admissible}
if $\varphi$ is real-valued and even.
We introduce the following norms on admissible functions in $\cS(\R)$
\begin{equation}\label{norm-varphi}
\cN_N(\varphi):= \sup_{u\ge 0} (1+u)^N \max_{0\le m\le N}|\varphi^{(m)}(u)|, \quad N\ge 0.
\end{equation}
\end{definition}

Observe that in the above we only need the values $\varphi(u)$
for $u \ge 0$.
Therefore, the condition ``$\varphi$ is even" can be replaced by
$\varphi^{(2\nu+1)}(0) =0$ for $\nu=0, 1, \dots$,
which implies that the even extension of $\varphi$ from $\R_+$ to $\R$ will have the required properties.


\begin{definition}\label{def:max-op-1}
Let $\varphi$ be an admissible function in $\cS(\R)$.
For any $f\in \cS'$ we define
\begin{equation}\label{def-max-1}
M(f;\varphi)(x):= \sup_{t>0} |\varphi(t\sqrt L)f(x)|,
\end{equation}
\begin{equation}\label{def-max-2}
M^*_a(f;\varphi)(x):= \sup_{t>0}\sup_{y\in M, \rho(x, y)\le at} |\varphi(t\sqrt L)f(y)|, \quad a\ge 1,
\end{equation}
and
\begin{equation}\label{def-max-3}
M^{**}_\gamma(f;\varphi)(x):= \sup_{t>0}\sup_{y\in M}
\frac{|\varphi(t\sqrt L)f(y)|}{\Big(1+\frac{\rho(x, y)}{t}\Big)^\gamma},\quad \gamma >0.
\end{equation}
\end{definition}

Observe that
\begin{equation}\label{relationship}
M(f;\varphi) \le M^*_a(f;\varphi) \le (1+a)^\gamma M^{**}_\gamma(f;\varphi), \quad \forall f\in\cS'.
\end{equation}

We now introduce the grand maximal operator.


\begin{definition}\label{def:max-op-2}
Denote
$$
\cF_N:= \{\varphi\in \cS(\R): \varphi \;{\rm is \; admissible \; and } \;\;\cN_N(\varphi)\le 1\}.
$$
The grand maximal operator is defined by
\begin{equation}\label{def-grand-max1}
\cM_N(f)(x):= \sup_{\varphi\in\cF_N} M^*_1(f;\varphi)(x), \quad f\in \cS',
\end{equation}
that is,
\begin{equation}\label{def-grand-max3}
\cM_{N}(f)(x):= \sup_{\varphi\in\cF_N}\sup_{t>0}\sup_{y\in M, \rho(x, y)\le t} |\varphi(t\sqrt L)f(y)|,
\end{equation}
where $N>0$ is sufficiently large $($to be specified$)$.
\end{definition}

It is readily seen that for any admissible function $\varphi$ and $a\ge 1$ one has
\begin{equation}\label{est-grand-max}
M^*_a(f;\varphi) \le a^N\cN_N(\varphi)\cM_{N}(f),
\quad \forall f\in\cS'.
\end{equation}

We will also use the following version of the Hardy-Littlewood maximal operator:
\begin{equation}\label{def-M}
M_\theta f(x):= \sup_{B\ni x} \Big(\frac{1}{|B|}\int_B|f(y)|^\theta d\mu(y)\Big)^{1/\theta},
\quad \theta >0.
\end{equation}

In the following we exhibit some important relations between the maximal operators.


\begin{proposition}\label{prop:max-operators}
Let $\varphi\in \cS(\R)$ be admissible and $\varphi(0)\ne 0$.

$(a)$ If $0< \theta \le 1$  and $\gamma > 2d/\theta$, then
\begin{equation}\label{int-MM1}
M^{**}_\gamma(f;\varphi)(x)
\le cM_\theta (M(f;\varphi))(x), \quad \forall f\in \cS',
\end{equation}
where $c= c(\theta, \gamma, d, \varphi)$.

$(b)$
If $0<\theta\le 1$ and $N > 6d/\theta+3d/2+2$, then
\begin{equation}\label{int-GMM1}
\cM_{N}(f)(x) \le cM_\theta (M(f;\varphi))(x), \quad \forall f\in \cS',
\end{equation}
where $c= c(\theta, d, \varphi)$.
\end{proposition}

For the proof of this proposition we need the following


\begin{lemma}\label{lem:Rych}
Suppose $\varphi\in \cS(\R)$ is admissible and $\varphi(0)=1$,
and let $N\ge 0$.
Then there exist even real-valued functions $\psi_0, \psi\in \cS(\R)$ such that
$\psi_0(0)=1$,
$\psi^{(\nu)}(0)=0$ for $\nu =0, 1, \dots, N$, and
for any $f\in \cS'$ and $j\in \bZ$
\begin{equation}\label{dec-Rych}
f = \psi_0(2^{-j}\sqrt L)\varphi(2^{-j}\sqrt L)f
+ \sum_{k=j}^\infty \psi(2^{-k}\sqrt L)[\varphi(2^{-k}\sqrt L) - \varphi(2^{-k+1}\sqrt L)]f,
\end{equation}
where the convergence is in $\cS'$.
\end{lemma}

\noindent
{\bf Proof.}
We borrow the idea for this proof from \cite[Theorem 1.6]{Rych}.
Evidently,
$$
\varphi(\lambda)^2 + \sum_{k=1}^\infty [\varphi(2^{-k}\lambda)^2- \varphi(2^{-k+1}\lambda)^2] =1,
\quad \lambda\in\R,
$$
and as $\varphi\in \cS(\R)$ the series converges absolutely.
From above
$$
1=\Big(\varphi(\lambda)^2
+ \sum_{k=1}^\infty \big[\varphi(2^{-k}\lambda)^2- \varphi(2^{-k+1}\lambda)^2\big]\Big)^N.
$$
It is easy to see that for $N\ge 1$ this identity can be written in the form
\begin{align*}
1& = \sum_{m=1}^N \binom{N}{m}\varphi(\lambda)^{2m}\big(1-\varphi(\lambda)^2\big)^{N-m}\\
&+ \sum_{k=1}^\infty
\sum_{m=1}^N \binom{N}{m}\big[\varphi(2^{-k}\lambda)^2- \varphi(2^{-k+1}\lambda)^2\big]^{m}
\big(1-\varphi(2^{-k}\lambda)^2\big)^{N-m},
\end{align*}
which leads to
\begin{equation}\label{dec-Rych-2}
\psi_0(\lambda)\varphi(\lambda)
+ \sum_{k=1}^\infty \psi(2^{-k}\lambda)[\varphi(2^{-k}\lambda) - \varphi(2^{-k+1}\lambda)] =1
\end{equation}
with
$$
\psi_0(\lambda):=
\sum_{m=1}^N \binom{N}{m}\varphi(\lambda)^{2m-1}\big(1-\varphi(\lambda)^2\big)^{N-m}
$$
and
$$
\psi(\lambda):= [\varphi(\lambda)+ \varphi(2\lambda)]
\sum_{m=1}^N \binom{N}{m}\big[\varphi(\lambda)^2- \varphi(2\lambda)^2\big]^{m-1}
\big(1-\varphi(\lambda)^2\big)^{N-m}.
$$
Clearly, $\psi_0, \psi\in \cS(\R)$, $\psi_0, \psi$ are even,
$\psi^{(\nu)}(0)=0$ for $\nu =0, 1, \dots, N-2$, and $\psi_0(0)=1$.

By replacing $\lambda$ in (\ref{dec-Rych-2}) by $\lambda/2$ and
subtracting the resulting identity from (\ref{dec-Rych-2}) we obtain
$$
\psi(2^{-1}\lambda)[\varphi(2^{-1}\lambda)-\varphi(\lambda)]
= \psi_0(2^{-1}\lambda)\varphi(2^{-1}\lambda)-\psi_0(\lambda)\varphi(\lambda).
$$
Hence, for any $f\in\cS'$ and $m>j$
\begin{align*}
&\psi_0(2^{-j}\sqrt L)\varphi(2^{-j}\sqrt L)f
+ \sum_{k=j}^m \psi(2^{-k}\sqrt L)[\varphi(2^{-k}\sqrt L) - \varphi(2^{-k+1}\sqrt L)]f\\
&\qquad\qquad
=\psi_0(2^{-m}\sqrt L)\varphi(2^{-m}\sqrt L)f \to f
\quad\hbox{as} \quad m\to \infty \quad \hbox{in $\;\cS'$},
\end{align*}
which implies (\ref{dec-Rych}).
Here we used Proposition~\ref{prop:convergence}.

Finally, by replacing $N$ with $N+2$ in the above proof we get what we need.
$\qed$

\medskip

\noindent
{\bf Proof of Proposition~\ref{prop:max-operators}.}
(a)
We borrow the idea for this proof from \cite[Lemma 3.2]{NS}.
Assume $0<\theta \le 1$ and $\gamma>2d/\theta$,
and let $f\in\cS'$.
We may assume that $\varphi(0)=1$ for otherwise we use $\varphi(0)^{-1}\varphi$ instead.
By Lemma~\ref{lem:Rych} there exist even real-valued functions $\psi_0, \psi\in \cS(\R)$ such that
$\psi_0(0)=1$,
$\psi^{(\nu)}(0)=0$ for $\nu =0, 1, \dots, N$, and for any $j\in \bZ$
(\ref{dec-Rych}) holds.

Fix $t>0$
and let $2^{-j} \le t< 2^{-j+1}$.
Using (\ref{dec-Rych}) we get
\begin{align*}
&\frac{|\varphi(t\sqrt L)f(y)|}{\big(1+\frac{\rho(x, y)}{t}\big)^\gamma}
\le c\frac{|\varphi(t\sqrt L)\psi_0(2^{-j}\sqrt L)\varphi(2^{-j}\sqrt L)f(y)|}{(1+2^j\rho(x, y))^\gamma}\\
& \qquad\qquad\qquad + c\sum_{k=j}^\infty
\frac{|\varphi(t\sqrt L)\psi(2^{-k}\sqrt L)[\varphi(2^{-k}\sqrt L)- \varphi(2^{-k+1}\sqrt L)]f(y)|}
{(1+2^j\rho(x, y))^\gamma}.
\end{align*}
Let $\omega(\lambda):= \varphi(t2^j\lambda)\psi(2^{-(k-j)}\lambda)$.
Then
$
\omega(2^{-j}\sqrt L)= \varphi(t\sqrt L)\psi(2^{-k}\sqrt L).
$

Now, choose $N > 3\gamma+3d/2+2$ and set $m:= \lfloor\gamma+d/2+1\rfloor$.
As $\varphi, \psi\in \cS$ there exists a constant $c>0$ such that
\begin{equation}\label{local-varphi}
|\varphi^{(\nu)}(\lambda)| \le c(1+\lambda)^{-N},
\quad
|\psi^{(\nu)}(\lambda)| \le c(1+\lambda)^{-N},
\quad \lambda>0, \;\;\nu=0, 1, \dots, N,
\end{equation}
yielding
$$ 
|\omega^{(\nu)}(\lambda)| \le c(1+\lambda)^{-N},
\quad \lambda>0, \;\;\nu=0, 1, \dots, N.
$$ 
From this estimate we obtain for $\lambda \ge 2^{(k-j)/2}$
$$
|\omega^{(\nu)}(\lambda)| \le c(1+\lambda)^{-m-d-1}2^{-(k-j)(N-m-d-1)/2}
$$
and using that $N \ge 3\gamma+3d/2+2 +2\eps$ for some $\eps>0$ it follows that
\begin{equation}\label{local-omega}
|\omega^{(\nu)}(\lambda)|
\le c2^{-(k-j)(\gamma+\eps)}(1+\lambda)^{-m-d-1},\;\; \lambda \ge 2^{(k-j)/2},
\;\;\nu=0, 1, \dots, N.
\end{equation}
On the other hand, as $\psi^{(\nu)}(0)=0$ for $\nu =0, 1, \dots, N$,
we use Taylor's formula and (\ref{local-varphi}) to obtain
$|\psi^{(\nu)}(\lambda)| \le c\lambda^{N-\nu}$, $\lambda>0$, $\nu=0, 1, \dots, N$.
Hence,
$$
\Big|\Big(\frac{d}{d\lambda}\Big)^\nu \psi(2^{-(k-j)}\lambda)\Big|
\le c2^{-(k-j)N}\lambda^{N-\nu}
\le c2^{-(k-j)N/2}
\quad\hbox{for} \;\;0\le \lambda \le 2^{(k-j)/2}.
$$
From this estimate and (\ref{local-varphi}) we infer 
$$ 
|\omega^{(\nu)}(\lambda)| \le c2^{-(k-j)N/2}(1+\lambda)^{-N},
\quad 0\le \lambda \le 2^{(k-j)/2},\;\; \nu=0, 1, \dots, N.
$$ 
In turn, this estimate and (\ref{local-omega}) imply that (\ref{local-omega})
holds for $0<\lambda <\infty$.
Now, Theorem~\ref{thm:S-local-kernels}, applied to $\omega(2^{-j}\sqrt L)$,
leads to the following estimate on the kernel of the operator $\varphi(t\sqrt L)\psi(2^{-k}\sqrt L)$
(recall that $2^{-j}\le t<2^{-j+1}$)
\begin{equation}\label{local-kernel}
|\varphi(t\sqrt L)\psi(2^{-k}\sqrt L)(x, y)|
\le \frac{c 2^{-(k-j)(\gamma+\eps)}}{|B(y, 2^{-j})|\big(1+2^j\rho(x, y)\big)^\gamma},
\quad x, y\in M.
\end{equation}
Consequently,
\begin{align*}
&\frac{|\varphi(t\sqrt L)\psi(2^{-k}\sqrt L)[\varphi(2^{-k}\sqrt L)- \varphi(2^{-k+1}\sqrt L)]f(y)|}
{(1+2^j\rho(x, y))^\gamma}\\
&\le c\int_M \frac{2^{-(\gamma+\eps)(k-j)}[|\varphi(2^{-k}\sqrt L)f(z)|+ |\varphi(2^{-k+1}\sqrt L)f(z)|]d\mu(z)}
{|B(z, 2^{-j})|\big(1+2^j\rho(y, z)\big)^\gamma\big(1+2^j\rho(x, y)\big)^\gamma}\\
&\le c\int_M \frac{2^{-(\gamma+\eps)(k-j)}[|\varphi(2^{-k}\sqrt L)f(z)|+ |\varphi(2^{-k+1}\sqrt L)f(z)|]d\mu(z)}
{|B(z, 2^{-j})|\big(1+2^j\rho(x, z)\big)^\gamma}\\
&\le c2^{-(k-j)\eps}\int_M \frac{[|\varphi(2^{-k}\sqrt L)f(z)|+ |\varphi(2^{-k+1}\sqrt L)f(z)|]d\mu(z)}
{|B(z, 2^{-j})|\big(1+2^k\rho(x, z)\big)^\gamma}\\
& \le c2^{-(k-j)\eps}\big[M^{**}_\gamma(f;\varphi)(x)\big]^{1-\theta}
\int_M \frac{[|\varphi(2^{-k}\sqrt L)f(z)|^\theta
+ |\varphi(2^{-k+1}\sqrt L)f(z)|^\theta]d\mu(z)}
{|B(z, 2^{-j})|\big(1+2^j\rho(x, z)\big)^{\gamma\theta}}.
\end{align*}
Similarly
\begin{align*}
&\frac{|\varphi(t\sqrt L)\psi_0(2^{-j}\sqrt L)\varphi(2^{-j}\sqrt L)f(y)|}{(1+2^j\rho(x, y))^\gamma}\\
&\qquad\qquad\qquad\qquad\qquad
\le c\big[M^{**}_\gamma(f;\varphi)(x)\big]^{1-\theta}
\int_M \frac{|\varphi(2^{-j}\sqrt L)f(z)|^\theta d\mu(z)}
{|B(z, 2^{-j})|\big(1+2^j\rho(x, z)\big)^{\gamma\theta}}.
\end{align*}
Putting the above estimates together we get
\begin{align*}
&\frac{|\varphi(t\sqrt L)f(y)|}{\big(1+\frac{\rho(x, y)}{t}\big)^\gamma}
\le c \big[M^{**}_\gamma(f;\varphi)(x)\big]^{1-\theta}\sum_{k=j}^\infty
2^{-(k-j)\eps}\int_M \frac{|\varphi(2^{-k}\sqrt L)f(z)|^\theta d\mu(z)}
{|B(z, 2^{-j})|\big(1+2^j\rho(x, z)\big)^{\gamma\theta}}.
\end{align*}
Using also (\ref{D2}), this yields
\begin{align*}
\big[M^{**}_\gamma(f;\varphi)(x)\big]^{\theta}
\le c \sum_{k=j}^\infty
2^{-(k-j)\eps}\int_M \frac{|\varphi(2^{-k}\sqrt L)f(z)|^\theta d\mu(z)}
{|B(x, 2^{-j})|\big(1+2^j\rho(x, z)\big)^{\gamma\theta-d}}.
\end{align*}
Denote briefly $F(z):= \varphi(2^{-k}\sqrt L)f(z)$.
We have
\begin{align*}
&\int_M \frac{|F(z)|^\theta d\mu(z)}
{|B(x, 2^{-j})|\big(1+2^j\rho(x, z)\big)^{\gamma\theta-d}}
= \int_{B(x, 2^{-j})} \frac{|F(z)|^\theta d\mu(z)}
{|B(x, 2^{-j})|\big(1+2^j\rho(x, z)\big)^{\gamma\theta-d}}\\
& \qquad\qquad\qquad
+ \sum_{m=1}^\infty
\int_{B(x, 2^{m-j})\setminus B(x, 2^{m-j-1})} \frac{|F(z)|^\theta d\mu(z)}
{|B(x, 2^{-j})|\big(1+2^j\rho(x, z)\big)^{\gamma\theta-d}}\\
&\qquad\qquad\qquad
\le c\sum_{m=0}^\infty
\frac{|B(x, 2^{m-j})|}{2^{m(\gamma\theta-d)}|B(x, 2^{-j})|}
\frac{1}{|B(x, 2^{m-j})|}\int_{B(x, 2^{m-j})}|F(z)|^\theta d\mu(z)\\
&\qquad\qquad\qquad
\le cM_\theta(F)(x)^\theta \sum_{m=0}^\infty
\frac{2^{md}}{2^{m(\gamma\theta-d)}}
\le cM_\theta(F)(x)^\theta.
\end{align*}
Here we used (\ref{doubling}) and that $\gamma > 2d/\theta$.
From above it follows that
\begin{align*}
\big[M^{**}_\gamma(f;\varphi)(x)\big]^{\theta}
\le c\sum_{k=j}^\infty 2^{-(k-j)\eps}M_\theta(\varphi(2^{-k}\sqrt L)f)(x)^\theta
\le cM_\theta (M(f;\varphi))(x)^\theta,
\end{align*}
which yields (\ref{int-MM1}).

\smallskip


To prove (b) we will proceed quite as in the proof of (a).
Let $\phi\in \cF_N$ and assume that $\varphi\in \cS(\R)$ is admissible.
Choose $\gamma > 2d/\theta$ so that $N > 3\gamma +3d/2+2$.
Then there exists $\eps>0$ such that $N \ge 3\gamma +3d/2+2+2\eps$.

Assume $t>0$
and let $2^{-j} \le t< 2^{-j+1}$.
Just as in the proof of (a),
by Lemma~\ref{lem:Rych} there exist even real-valued functions $\psi_0, \psi\in \cS(\R)$ such that
$\psi^{(\nu)}(0)=0$ for $\nu =0, 1, \dots, N$ and for any $j\in \bZ$
(\ref{dec-Rych}) holds.
%
Hence, for $f\in \cS'$,
\begin{align*}
&\frac{|\phi(t\sqrt L)f(y)|}{\big(1+\frac{\rho(x, y)}{t}\big)^\gamma}
\le c\frac{|\phi(t\sqrt L)\psi_0(2^{-j}\sqrt L)\varphi(2^{-j}\sqrt L)f(y)|}{(1+2^j\rho(x, y))^\gamma}\\
& \qquad\qquad\qquad + c\sum_{k=j}^\infty
\frac{|\phi(t\sqrt L)\psi(2^{-k}\sqrt L)[\varphi(2^{-k}\sqrt L)- \varphi(2^{-k+1}\sqrt L)]f(y)|}
{(1+2^j\rho(x, y))^\gamma}.
\end{align*}
Just as in (\ref{local-kernel}) we have
$$
|\phi(t\sqrt L)\psi(2^{-k}\sqrt L)(x, y)|
\le \frac{c 2^{-(k-j)(\gamma+\eps)}}{|B(y, 2^{-j})|\big(1+2^j\rho(x, y)\big)^\gamma},
$$
where the constant $c>0$ is independent of $\phi$ due to
$\cN_N(\phi) \le 1$.
Therefore, as in the proof of (a)
\begin{align*}
&\frac{|\phi(t\sqrt L)\psi(2^{-k}\sqrt L)[\varphi(2^{-k}\sqrt L)- \varphi(2^{-k+1}\sqrt L)]f(y)|}
{(1+2^j\rho(x, y))^\gamma}\\
&\le c\int_M \frac{2^{-(\gamma+\eps)(k-j)}[|\varphi(2^{-k}\sqrt L)f(z)|+ |\varphi(2^{-k+1}\sqrt L)f(z)|]d\mu(z)}
{|B(z, 2^{-j})|\big(1+2^j\rho(y, z)\big)^\gamma\big(1+2^j\rho(x, y)\big)^\gamma}\\
%
%
& \le c2^{-(k-j)\eps}\big[M^{**}_\gamma(f;\varphi)(x)\big]^{1-\theta}
\int_M \frac{[|\varphi(2^{-k}\sqrt L)f(z)|^\theta
+ |\varphi(2^{-k+1}\sqrt L)f(z)|^\theta]d\mu(z)}
{|B(z, 2^{-j})|\big(1+2^j\rho(x, z)\big)^{\gamma\theta}}.
\end{align*}
Similarly
\begin{align*}
&\frac{|\phi(t\sqrt L)\psi_0(2^{-j}\sqrt L)\varphi(2^{-j}\sqrt L)f(y)|}{(1+2^j\rho(x, y))^\gamma}\\
&\qquad\qquad\qquad\qquad\qquad
\le c\big[M^{**}_\gamma(f;\varphi)(x)\big]^{1-\theta}
\int_M \frac{|\varphi(2^{-j}\sqrt L)f(z)|^\theta d\mu(z)}
{|B(z, 2^{-j})|\big(1+2^j\rho(x, z)\big)^{\gamma\theta}}.
\end{align*}
Here the constant $c>0$ is independent of $\phi$ since
$\cN_N(\phi) \le 1$.
As before denoting $F(z):= \varphi(2^{-k}\sqrt L)f(z)$ we obtain
\begin{align*}
\int_M \frac{|F(z)|^\theta d\mu(z)}
{|B(z, 2^{-j})|\big(1+2^j\rho(x, z)\big)^{\gamma\theta}}
\le c\int_M \frac{|F(z)|^\theta d\mu(z)}
{|B(x, 2^{-j})|\big(1+2^j\rho(x, z)\big)^{\gamma\theta-d}}
\le cM_\theta(F)(x)^\theta.
\end{align*}
Therefore,
\begin{align*}
\frac{|\phi(t\sqrt L)f(y)|}{\big(1+\frac{\rho(x, y)}{t}\big)^\gamma}
&\le c \big[M^{**}_\gamma(f;\varphi)(x)\big]^{1-\theta}\sum_{k=j}^\infty
2^{-(k-j)\eps}\int_M \frac{|\varphi(2^{-k}\sqrt L)f(z)|^\theta d\mu(z)}
{|B(z, 2^{-j})|\big(1+2^j\rho(x, z)\big)^{\gamma\theta}}\\
&\le c \big[M^{**}_\gamma(f;\varphi)(x)\big]^{1-\theta}\sum_{k=j}^\infty
2^{-(k-j)\eps}\big[M_\theta (|\varphi(2^{-k}\sqrt L)f|)(x)\big]^\theta\\
& \le c \big[M^{**}_\gamma(f;\varphi)(x)\big]^{1-\theta}
\big[M_\theta (M(f;\varphi))(x)\big]^\theta\sum_{k=j}^\infty 2^{-(k-j)\eps}\\
&\le cM_\theta (M(f;\varphi))(x),
\end{align*}
where for the last estimate we used that
$M^{**}_\gamma(f;\varphi)(x) \le cM_\theta (M(f;\varphi))(x)$, by (\ref{int-GMM1}).
Thus
\begin{align*}
\sup_{t>0}\sup_{y\in M, \rho(x, y)\le t} |\phi(t\sqrt L)f(y)|
\le 2^\gamma\sup_{t>0}\sup_{y\in M}\frac{|\phi(t\sqrt L)f(y)|}{\big(1+\frac{\rho(x, y)}{t}\big)^\gamma}
\le cM_\theta (M(f;\varphi))(x),
\end{align*}
which completes the proof.
$\qed$


\begin{definition}\label{def:Hp}
The Hardy space $H^p$, $0<p\le 1$, is defined as the set of all
distributions $f\in \cS'$ such that
\begin{equation}\label{def-Hp}
\|f\|_{H^p}:= \big\|\sup_{t>0}|e^{-t^2L}f(\cdot)|\big\|_{L^p} <\infty.
\end{equation}
\end{definition}

Proposition~\ref{prop:max-operators} leads to the following


\begin{theorem}\label{thm:Hp}
Let $0<p\le 1$.
Then for any $N >6d/p+3d/2+2$, $\gamma >2d/p$, $a\ge 1$, and
an admissible $\varphi\in \cS(\R)$ with $\varphi(0)\ne 0$ we have for all $f\in \cS'$
\begin{equation}\label{equiv-Hp-norms}
\|f\|_{H^p} \sim \|\cM_{N}(f)\|_{L^p}
\sim \|M(f;\varphi)\|_{L^p}
\sim \|M^*_a(f;\varphi)\|_{L^p}
\sim \|M^{**}_\gamma(f;\varphi)\|_{L^p}.
\end{equation}
Here the constants in the equivalences involving $\varphi$ depend
not only on the parameters but on $\varphi$ as well.
\end{theorem}

\noindent
{\bf Proof.}
Write $\Phi(\lambda):= e^{-\lambda^2}$. Apparently $\Phi\in\cS(\R)$, $\Phi$ is admissible, and $\Phi(0)\ne 0$.
Let $N> 6d/p+3d/2+2$ and choose $\theta$ so that $0<\theta <p$ and $N> 6d/\theta+3d/2+2$.
Then applying Proposition~\ref{prop:max-operators} (b) we get
$$
\|\cM_{N}(f)\|_{p} \le c \|M_\theta (M(f;\Phi))\|_p \le c \|M(f;\Phi)\|_p =c \|f\|_{H^p},
$$
where we used the maximal inequality:
$ \|M_\theta f\|_p \le c \|f\|_p$ if $0<\theta<p$, see \cite{Stein}.\\
In the other direction, using (\ref{relationship}) and (\ref{est-grand-max}) we get
$$
\|f\|_{H^p}=\|M(f;\Phi)\|_p \le \|M_1^*(f;\Phi)\|_p \le c\|\cM_{N}(f)\|_{p}.
$$
Thus the first equivalence in (\ref{equiv-Hp-norms}) is established.

Just in the same way we get
$\|\cM_{N}(f)\|_{L^p}\sim \|M(f;\varphi)\|_{L^p}$
with constants of equivalence depending in addition on $\varphi$.
We choose $\theta$ so that $0<\theta<p$ and $\gamma > 2d/\theta$
and apply Proposition~\ref{prop:max-operators} (a) and the maximal inequality as above
to obtain
$\|M^{**}_\gamma(f;\varphi)\|_{L^p} \le c \|M(f;\varphi)\|_{L^p}$.
All other estimates we need follow from (\ref{relationship}) and (\ref{est-grand-max}).
$\qed$

\subsection{Characterization of \boldmath $H^p$ via the Poisson semigroup}

We next establish a characterization of the Hardy spaces $H^p$ 
with the r\^{o}le of $e^{-t^2L}$ in (\ref{def-Hp}) replaced by $e^{-t\sqrt{L}}$.
To this end we first show that $H^p$ consists of bounded distributions.


\begin{definition}\label{def:bounded}
We say that $f\in\cS'$ is a bounded distribution if there exists $N\ge 0$ such that
for every admissible $\varphi\in\cS(\R)$ we have
$\varphi(\sqrt{L})f \in L^\infty(M)$
and $\|\varphi(\sqrt{L})f\|_\infty \le c(f)\cN_N(\varphi)$,
where $c(f)>0$ depends only on $f$; $\cN_N(\varphi)$ is from~$(\ref{norm-varphi})$.
\end{definition}


\begin{lemma}\label{lem:bounded}
If $f\in H^p$, $0<p\le 1$, then $f$ is a bounded distribution
and for every admissible function $\varphi\in\cS(\R)$ we have
$\|\varphi(\sqrt{L})f\|_\infty \le c\cN_N(\varphi)\|f\|_{H^p}$
with $N$ as in Theorem~\ref{thm:Hp}.
\end{lemma}

\noindent
{\bf Proof.}
Let $f\in H^p$, $0<p\le 1$.
Then, by Theorem~\ref{thm:Hp}, we have  $\|\cM_{N}(f)\|_{L^p}\sim \|f\|_{H^p}$ provided $N>6d/p+3d/2+2$.
Let $\varphi\in\cS(\R)$ be admissible. 
Then by (\ref{est-grand-max}) we have
$M^*_1(f; \varphi) \le \cN_{N} (\varphi)\cM_{N}(f)$.
Therefore, for each $x\in M$
\begin{align*}
|\varphi(\sqrt{L})f(x)|^p
&\le \inf_{y: \rho(x, y)\le 1}\sup_{z: \rho(z, y)\le 1}|\varphi(\sqrt{L})f(z)|^p
\le \inf_{y: \rho(x, y)\le 1} M^*_1(f; \varphi)(y)^p\\
&\le \frac{1}{|B(x, 1)|}\int_{B(x, 1)}M^*_1(f; \varphi)^p(y) d\mu(y)\\
&\le c\cN_{N}(\varphi)^p\int_{M}\cM_{N}(f)(y)^p d\mu(y)
\le c\cN_{N}(\varphi)^p\|f\|_{H^p}^p,
\end{align*}
where we used 
the non-collapsing condition (\ref{non-collapsing}).
The lemma follows.
$\qed$

The above lemma enables us to identify $e^{-t\sqrt{L}}f$ for $f\in H^p$
with a well-defined bounded function.


\begin{definition}\label{def:Poisson-Hp}
Let $\varphi\in C^\infty(\R)$ be even and real-valued, $\varphi(\lambda)=1$ for $x\in [-1, 1]$
and $\supp \varphi \subset [-2, 2]$.
Set $\theta(\lambda):= e^{-\lambda}(1-\varphi(\lambda))$ for $\lambda \ge 0$
and denote again by $\theta$ the even extension of $\theta$ on $\R$.
Clearly, $\theta\in \cS(\R)$ is admissible.
Given a bounded distribution $f\in\cS'$,
for example
$f\in H^p$, $0<p\le 1$, we define
\begin{equation}\label{def-Poisson}
e^{-t\sqrt{L}}f := e^{-t\sqrt{L}}\varphi(t\sqrt{L})f + \theta(t\sqrt{L})f,
\quad t>0.
\end{equation}
\end{definition}



\begin{lemma}\label{lem:well-def}
If $f\in \cS'$ is a bounded distribution $($Definition~\ref{def:bounded}$)$,
then $e^{-t\sqrt{L}}f$ is a well-defined function and $e^{-t\sqrt{L}}f\in L^\infty(M)\cap C(M)$.
\end{lemma}

\noindent
{\bf Proof.}
We will use the subordination formula:
For any $f\in L^2(M)$
\begin{equation}\label{subordinat-1}
e^{-t\sqrt{L}}f
%
= \frac{1}{\sqrt{\pi}}\int_0^\infty\frac{t}{s^2}e^{-\frac{t^2}{4s^2}}e^{-s^2L}f ds, \quad t>0,
\end{equation}
which follows easily from the spectral $L^2$-theory.
Let $\varphi$ be as in Definition~\ref{def:Poisson-Hp}.
Clearly, $e^{-\lambda}\varphi(\lambda)$ is bounded and compactly supported and
by \cite[Theorem~3.7]{CKP} it follows that $e^{-t\sqrt{L}}\varphi(t\sqrt{L})$
is an integral operator with bounded and continuous kernel.
By Theorem~\ref{thm:S-local-kernels} it follows that the kernel of the operator $\theta(t\sqrt{L})$
is also bounded and continuous.
Hence, in light of (\ref{def-Poisson}), the kernel $e^{-t\sqrt L}(x, y)$ of the operator $e^{-t\sqrt L}$
is continuous and bounded.
Using also the fact that the heat kernel $e^{-tL}(x, y)=p_t(x, y)$ is H\"{o}lder continuous (see (\ref{lip}))
it readily follows from (\ref{subordinat-1}) that 
\begin{equation}\label{subordinat-2}
e^{-t\sqrt{L}}(x,y)
= \frac{1}{\sqrt{\pi}}\int_0^\infty\frac{t}{s^2}e^{-\frac{t^2}{4s^2}}e^{-s^2L}(x, y) ds,
\quad \forall x, y\in M.
\end{equation}
Now, taking into account the fact that for any $f\in L^\infty(M)$
$$
\int_0^\infty\int_M
|e^{-\frac{t^2}{4s^2}}e^{-s^2L}(x, y)|
|f(y)| d\mu(y)ds \le c\|f\|_\infty
$$
(\ref{subordinat-2}) yields that for any $f\in L^\infty(M)$
\begin{equation}\label{subordinat-3}
e^{-t\sqrt{L}}f(x)
= \frac{1}{\sqrt{\pi}}\int_0^\infty\frac{t}{s^2}e^{-\frac{t^2}{4s^2}}e^{-s^2L}f(x) ds, \quad x\in M,
\end{equation}
and
$\|e^{-t\sqrt{L}}f\|_\infty \le c\|f\|_\infty$.

Assume that $f\in\cS'$ is a bounded distribution and let $\varphi$ be as in Definition~\ref{def:Poisson-Hp}.
Then
$\|\varphi(\sqrt{L})f\|_\infty \le c(f)\cN_N(\varphi)$
for some $N\ge 0$.
From (\ref{norm-varphi}), it is easy to see that
$\cN_N(\varphi (t\cdot)) \le (t+t^{-1})^N \cN_N(\varphi)$, $t>0$,
and hence using the above
\begin{equation}\label{estim-1}
\|e^{-t\sqrt{L}}\varphi(t\sqrt{L})f\|_\infty \le c(f)(t+t^{-1})^N \cN_N(\varphi).
\end{equation}
On the other hand, the function $\theta$ from Definition~\ref{def:Poisson-Hp} is admissible in $\cS(\R)$
and hence $\|\theta(t\sqrt{L})f\|_\infty \le c(f)(t+t^{-1})^N\cN_N(\theta)$, $t>0$.
This coupled with (\ref{estim-1}) implies that
 $e^{-t\sqrt{L}}f\in L^\infty(M)$.


We next show that $e^{-t\sqrt{L}}f \in C(M)$.
Write $F:=\varphi(t\sqrt{L})f$. Then by (\ref{subordinat-3})
\begin{align*}
&e^{-t\sqrt{L}}F(x)-e^{-t\sqrt{L}}F(x')
= \frac{1}{\sqrt{\pi}}\int_0^\infty\frac{t}{s^2}e^{-\frac{t^2}{4s^2}}
\big[e^{-s^2L}F(x)-e^{-s^2L}F(x')\big]ds\\
&\qquad\qquad\qquad
= \frac{1}{\sqrt{\pi}}\int_0^\infty\frac{t}{s^2}e^{-\frac{t^2}{4s^2}}
\int_M\big[e^{-s^2L}(x, y)-e^{-s^2L}(x', y)\big]F(y)d\mu(y)ds.
\end{align*}
From (\ref{lip}) and (\ref{D2}), 
for any $\sigma>0$ there exists a constant $c_\sigma>0$ such that
\begin{align}\label{est-heat-1}
\big|e^{-s^2L}(x, y)-e^{-s^2L}(x', y)\big|
\le c_\sigma \big(\rho(x, x')/s\big)^\alpha
|B(x, s)|^{-1}\big(1+s^{-1}\rho(x,y)\big)^{-\sigma}
\end{align}
whenever $\rho(x, x') \le s$, and
\begin{align}\label{est-heat-2}
\big|e^{-s^2L}(x, y)\big|
\le c_\sigma |B(x,s)|^{-1}\big(1+s^{-1}\rho(x,y)\big)^{-\sigma}.
\end{align}
We choose $\sigma >2d$. Denote $A:= \rho(x, x')$ and assume $A>0$.
Then from above
\begin{align*}
|e^{-t\sqrt{L}}F(x)-e^{-t\sqrt{L}}F(x')|
&\le \int_0^A\cdots + \int_A^\infty \cdots
=: J_1+J_2.
\end{align*}
To estimate $J_1$ we use (\ref{est-heat-2}) and (\ref{tech-1}). We get
\begin{align*}
J_1 &\le c \|F\|_\infty
\int_0^A\frac{t}{s^2}e^{-\frac{t^2}{4s^2}}
\int_M\frac{d\mu(y)}{|B(x, s)|\big(1+s^{-1}\rho(x,y)\big)^\sigma}ds\\
&  
+ c \|F\|_\infty
\int_0^A\frac{t}{s^2}e^{-\frac{t^2}{4s^2}}
\int_M\frac{d\mu(y)}{|B(x', s)|\big(1+s^{-1}\rho(x',y)\big)^\sigma}ds\\
&  
\le c \|F\|_\infty \int_0^A\frac{t}{s^2}e^{-\frac{t^2}{4s^2}}ds
= c \|F\|_\infty \int_{t^2/4A}^\infty \frac{e^{-u}}{\sqrt u}du\\
& \le c \|F\|_\infty At^{-1}\int_0^\infty e^{-u}du
\le c \rho(x, x'),
\quad c=c(f, \varphi, t).
\end{align*}
Here for the equality we applied the substitution $u=t^2/4s^2$.

To estimate $J_2$ we use (\ref{est-heat-1}) and (\ref{tech-1}) and obtain
\begin{align*}
J_2 &\le c\rho(x, x')^\alpha \|F\|_\infty
\int_A^\infty\frac{t}{s^{2+\alpha}}e^{-\frac{t^2}{4s^2}}
\int_M\frac{d\mu(y)}{|B(x, s)|\big(1+s^{-1}\rho(x,y)\big)^\sigma}ds\\
&
\le c\rho(x, x')^\alpha \|F\|_\infty\int_0^\infty\frac{t}{s^{2+\alpha}}e^{-\frac{t^2}{4s^2}}ds
\le c \rho(x, x')^\alpha.
\qquad \hbox{$(u=t^2/4s^2)$}
\end{align*}
Therefore,
$$ 
|e^{-t\sqrt{L}}\varphi(t\sqrt{L})f(x)-e^{-t\sqrt{L}}\varphi(t\sqrt{L})f(x')|
\le c \rho(x, x')^\alpha, \quad \hbox{if}\;\; \rho(x, x') \le 1.
$$ 
Also, by Proposition~\ref{prop:dist} $\theta(t\sqrt L)f\in C(M)$ and hence $e^{-t\sqrt{L}}f \in C(M)$.

Furthermore, it is easy to see that the definition of $e^{-t\sqrt{L}}f$ in (\ref{def-Poisson})
is independent of the particular selection of the function $\varphi$.
$\qed$

\smallskip

We now come to the main point in this subsection.


\begin{theorem}\label{thm:poisson}
Let $f\in \cS'$ be a bounded distribution $($see Definition~\ref{def:bounded}$)$.
Then $f\in H^p$, $0<p\le 1$, if and only if
$\big\|\sup_{t>0}\big|e^{-t\sqrt{L}}f(\cdot)\big|\big\|_{L^p}<\infty$,
and
\begin{equation}\label{Hp-poisson}
\big\|\sup_{t>0}\big|e^{-t\sqrt{L}}f(\cdot)\big|\big\|_{L^p} \sim \|f\|_{H^p}.
\end{equation}
\end{theorem}

\noindent
{\bf Proof.}
Let $f\in H^p$, $0<p\le 1$.
We will show that the subordination formula (\ref{subordinat-1}) holds for this distribution pointwise:
\begin{equation}\label{subordinat-4}
e^{-t\sqrt{L}}f(x)
= \frac{1}{\sqrt{\pi}}\int_0^\infty\frac{t}{s^2}e^{-\frac{t^2}{4s^2}}e^{-s^2L}f(x) ds,
\quad\forall x\in M.
\end{equation}
Let $\psi \in C^\infty_0(\R)$ be just
as the function $\varphi$ from Definition~\ref{def:Poisson-Hp}.
By Lemma~\ref{lem:bounded} $f$ is a bounded distribution and hence, by lemma~\ref{lem:well-def},
$\psi(\delta\sqrt{L})f\in L^\infty(M)\cap C(M)$ for $\delta>0$,
taking into account Proposition~\ref{prop:dist} as well.
Then using (\ref{subordinat-3}) we infer
\begin{align}\label{subordinat-5}
e^{-t\sqrt{L}}\psi(\delta\sqrt{L})f(x)
&= \frac{1}{\sqrt{\pi}}\int_0^\infty\frac{t}{s^2}e^{-\frac{t^2}{4s^2}}e^{-s^2L}\psi(\delta\sqrt{L})f(x) ds,
\quad\forall x\in M.
\end{align}
However, by Proposition~\ref{prop:convergence} if $g\in L^\infty(M)\cap C(M)$, then
$\psi(\delta\sqrt{L})g(x) \to g(x)$, $\forall x\in M$
and by Lemma~\ref{lem:well-def} $e^{-t\sqrt{L}}f \in L^\infty(M)\cap C(M)$.
Therefore,
$$
\lim_{\delta\to 0}e^{-t\sqrt{L}}\psi(\delta\sqrt{L})f(x)
= \lim_{\delta\to 0}\psi(\delta\sqrt{L})e^{-t\sqrt{L}}f(x) = e^{-t\sqrt{L}}f(x),
\quad \forall x\in M.
$$
Here, to justify the equality
$e^{-t\sqrt{L}}\psi(\delta\sqrt{L})f(x) = \psi(\delta\sqrt{L})e^{-t\sqrt{L}}f(x)$
we use the fact that
$\langle e^{-t\sqrt{L}} g, h\rangle = \langle g, e^{-t\sqrt{L}} h\rangle$
for all $g\in L^\infty$, $h\in L^1$, which follows from (\ref{subordinat-2});
we also use (\ref{comute}), (\ref{def-Poisson}), and the fact that $f$ is a bounded distribution.

Similarly,
$
\lim_{\delta\to 0} e^{-s^2L}\psi(\delta\sqrt{L})f(x)
= e^{-s^2L}f(x),
$
$\forall x\in M$,
and in addition
$$
\|e^{-s^2L}\psi(\delta\sqrt{L})f\|_\infty
\le c\|\psi(\delta\sqrt{L})f\|_\infty \le c\|f\|_{H^p}.
$$
Passing to the limit in (\ref{subordinat-5}) as $\delta \to 0$, we obtain (\ref{subordinat-4});
we apply the dominated convergence theorem to justify the convergence of the integral in (\ref{subordinat-5})
as $\delta\to 0$.

Estimate (\ref{subordinat-4}) implies
$
\sup_{t>0}\big|e^{-t\sqrt{L}}f(x)\big| \le c \sup_{t>0}\big|e^{-t^2L}f(x)\big|
$
and hence
$$
\big\|\sup_{t>0}\big|e^{-t\sqrt{L}}f(\cdot)\big|\big\|_{L^p}
\le c \big\|\sup_{t>0}\big|e^{-t^2L}f(\cdot)\big|\big\|_{L^p}
= c\|f\|_{H^p}.
$$
This completes the first part of the proof.


\smallskip

For the other direction, we borrow the following claim from
\cite{Stein-1970}, pp. 182-183, see also \cite{Stein}, p. 99:
 There exists a real-valued function $\eta\in C[1, \infty)$
such that
$\eta(s)\le O(s^{-\gamma})$ as $s\to \infty$, $\forall \gamma>0$,
\begin{equation}\label{def-eta}
\int_1^\infty s^k\eta(s)ds = 0 \quad{ for }\;\; k=1, 2, \dots,
\quad\hbox{and}\;\; \int_1^\infty \eta(s)ds =1.
\end{equation}
The function
$\eta(s):= \frac{e}{\pi s} {\rm Im} \{\exp [-\omega(s-1)^{1/4}]\}$ with
$\omega=e^{-\pi i/4}$ will do \cite{Stein-1970}.

Consider the function
$\Phi(\lambda):= \int_1^\infty \eta(s) e^{-s|\lambda|}ds$ for $\lambda \in \R$.
It is easy to see that $\Phi\in \cS(\R)$, $\Phi$ is even, and $\Phi(0)=1$.
Therefore, from the $L^2$-theory, for any $f\in L^2$
\begin{equation}\label{Phi-eta}
\Phi(t\sqrt{L})f= \int_1^\infty \eta(s) e^{-st\sqrt{L}}fds.
\end{equation}

Assume that $f\in \cS'$ is a bounded distribution.
Then exactly as in the proof of the first part we use Lemma~\ref{lem:well-def} to show that
(\ref{Phi-eta}) is valid point-wise for this distribution, which 
implies
$\sup_{t>0}|\Phi(t\sqrt{L})f(x)| \le c\sup_{t>0}\big|e^{-t\sqrt{L}}f(x)\big|$.
Therefore,
$$
\|f\|_{H^p} \le c\big\|\sup_{t>0}|\Phi(t\sqrt{L})f(\cdot)|\big\|_{L^p}
\le c \big\|\sup_{t>0}\big|e^{-t\sqrt{L}}f(\cdot)\big|\big\|_{L^p}
$$
as desired.
The proof is complete.
$\qed$

\subsection{Some simple facts about \boldmath $H^p$ spaces}\label{subsec:facts}

Here we present without proof some simple facts about the Hardy spaces $H^p$, $0<p\le 1$,
in the setting of this paper.

(a) $H^p$ is continuously embedded in $\cS'$, that is,
for any $0<p\le 1$ there exist constants $m >0$ and $c>0$ such that
for every $f\in H^p$ and $\phi\in\cS$ one has
$|\langle f, \phi\rangle| \le c \cP_m(\phi)\|f\|_{H^p}$.

(b) $H^p$ is a complete quasi-normed space.

(c) For any $q >p$ the space $L^q(M)$ is dense in $H^p$,
moreover, $\cS$ is dense in $H^p$.

Observe also that as in the classical case on $\R^n$ one has
$H^p=L^p$ with equivalent norms whenever $p>1$.


\section{Atomic Hardy spaces}\label{sec:atomic-Hp-spaces}
\setcounter{equation}{0}

We consider two versions of atomic Hardy spaces depending on whether
$\mu(M) =\infty$ or $\mu(M) <\infty$.
In defining the atomic Hardy spaces we borrow from \cite{HLMMY, Duong-Li}.

\subsection{Atomic Hardy spaces in the noncompact case}\label{sec:Hom-atomic-Hp-spaces}


\begin{definition} \label{def-atoms-hom}
Let $0<p\le 1$ and $n:=\lfloor d/2p\rfloor +1$, where $d$ is from the doubling property $(\ref{doubling})$.
A~function $a(x)$ is called an atom associated with the operator $L$ if
there exists a function $b\in D(L^n)$ and a ball $B$ of radius $r=r_B >0$
such that
\begin{enumerate}
\item[(i)] $a=L^nb,$
\item[(ii)] $\supp L^k b\subset B$, $k=0, 1, \dots, n$, and
\item[(iii)] $\|L^k b\|_\infty \le r^{2(n-k)}|B|^{-1/p}$, $k=0, 1, \dots, n$.
\end{enumerate}

\end{definition}


\begin{definition} \label{def-Hardy-Atomic-hom}
The atomic Hardy space $H^p_A$, $0<p\le 1$, is defined as the set of all distributions $f\in \cS'$
that can be represented in the form
\begin{equation}\label{atomic-Hardy-hom}
f=\sum_{j=1}^\infty \lambda_j a_j,
\quad\hbox{where}\quad
\sum_{j=1}^\infty |\lambda_j|^p <\infty,
\end{equation}
$\{a_j\}$ are atoms, and the convergence is in $\cS'$.
We set
\begin{equation}\label{norm-atomic-Hardy-hom}
\|f\|_{H^p_A}:=\inf_{f=\sum_{j\ge 1} \lambda_j a_j} \Big(\sum_{j\ge 1} |\lambda_j|^p\Big)^{1/p},
\quad f\in H^p_A.
\end{equation}
\end{definition}

\subsection{Atomic Hardy spaces in the compact case}\label{sec:inhom-atomic-Hp-spaces}

In the case when $\mu(M)<\infty$, we use the atoms from above with the addition of one more kind of atoms,
say, $A \in L^\infty(M)$, with the property:
\begin{equation}\label{def-A}
\|A\|_\infty \le |M|^{-1/p}.
\end{equation}
Then the atomic Hardy space $H^p_A$, $0<p\le 1$, is defined just as in the noncompact case above.

\section{Equivalence of maximal and atomic Hardy spaces}\label{sec:equivalence-Hp-spaces}
\setcounter{equation}{0}

We now come to the main result of this article.


\begin{theorem}\label{thm:hardy}
In the setting of this paper,
we have
$H^p=H^p_A$, $0<p\le 1$, and
\begin{equation}\label{hardy}
\|f\|_{H^p_A} \sim \|f\|_{H^p}
\quad\hbox{for}\quad
f\in H^p.
\end{equation}

\end{theorem}

We will first carry out the proof of this theorem in the noncompact case
and then explain the modifications that need to me made in the compact case.

\subsection{Proof of the embedding \boldmath $H^p \subset H^p_A$ in the noncompact case}
\label{subsec:embed-Hp-HpA}

We next show that if $f\in H^p$, $0<p\le 1$, then $f\in H^p_A$ and
$\|f\|_{H^p_A} \le c\|f\|_{H^p}$.


We begin with a simple decomposition identity which will play a central r\^{o}le in this proof;
it relies on the following


\begin{lemma}\label{lem:construction}
For any $m\ge 1$ there exists a function
$\varphi\in \cS(\R)$ 
such that
$\varphi$ is real-valued and even,
$\supp \hat \varphi \subset [-1, 1]$,
$\varphi(0)=1$, and
$\varphi^{(\nu)}(0)=0$ for $\nu =1, 2, \dots, m$.
\end{lemma}

To make our exposition more fluid we relegate the proof of this lemma to the appendix.

\smallskip

Assume $\varphi\in \cS(\R)$, $\varphi$ is real-valued and even,
$\supp \hat \varphi \subset [-1, 1]$,
$\varphi(0)=1$, and
$\varphi^{(\nu)}(0)=0$ for $\nu =1, 2, \dots, \NN-1$,
where $\NN$ is even and sufficiently large.
More precisely, we choose $\NN \ge 2n+d+d/p+1$ with $n$ from Definition \ref{def-atoms-hom}.
The existence of such $\varphi$ is guaranteed by Lemma~\ref{lem:construction}.
Appealing to Proposition~\ref{prop:convergence},
for any $f\in\cS'$ we have $f=\lim_{j\to\infty}\varphi(2^{-j}\sqrt L)^2f$ in $\cS'$,
which implies this representation:
For any $j\in\bZ$
\begin{align*}
f&=\varphi(2^{-j}\sqrt L)^2f + \sum_{k=j}^\infty \big[\varphi(2^{-k-1}\sqrt L)^2f-\varphi(2^{-k}\sqrt L)^2f\big]\\
&=\varphi(2^{-j}\sqrt L)^2f + \sum_{k=j}^\infty
\big[\varphi(2^{-k-1}\sqrt L)-\varphi(2^{-k}\sqrt L)\big]
\big[\varphi(2^{-k-1}\sqrt L)+\varphi(2^{-k}\sqrt L)\big]f.
\end{align*}
Setting $\psi(\lambda) := \varphi(\lambda)-\varphi(2\lambda)$ and
$\tilde\psi(\lambda) := \varphi(\lambda)+\varphi(2\lambda)$,
we arrive at
\begin{equation}\label{main-repr}
f= \varphi(2^{-j}\sqrt L)^2f + \sum_{k=j+1}^\infty \psi(2^{-k}\sqrt L)\tilde\psi(2^{-k}\sqrt L)f,
\quad f\in \cS'.
\end{equation}
Clearly, $\psi,\tilde\psi\in\cS(\R)$, $\psi,\tilde\psi$ are even,
$\supp\hat\psi \subset [-2, 2]$ and $\supp \hat{\tilde\psi} \subset [-2, 2]$.
Then by the final speed propagation property (Proposition~\ref{prop:finite-sp})
\begin{equation}\label{supp-psi}
\supp \psi(2^{-k}\sqrt L)(x, \cdot) \subset B(x, \cc2^{-k}),
\quad
\supp \tilde\psi(2^{-k}\sqrt L)(x, \cdot) \subset B(x, \cc2^{-k}),
\end{equation}
where $\cc >1$ is a constant.

For later use, observe also that
\begin{equation}\label{derivatives}
\psi^{(\nu)}(0) =0 \quad\hbox{for}\quad \nu=0, 1, \dots, \NN-1.
\end{equation}

From now on we will use the following more compact notation: 
\begin{equation}\label{brief-notation}
\varphi_k:=\varphi(2^{-k}\sqrt L), \quad
\psi_k:=\psi(2^{-k}\sqrt L),
\quad\hbox{and}\quad
\tilde\psi_k:=\tilde\psi(2^{-k}\sqrt L).
\end{equation}
The kernels of these operators will be denoted by
$\varphi_k(x, y)$, $\psi_k(x, y)$, and $\tilde\psi_k(x, y)$.
Observe that since $\varphi$, $\psi$, and $\tilde\psi$ are real-valued
we have $\varphi_k(y, x)=\varphi_k(x, y)$
and similarly for the others.
By Theorem~\ref{thm:S-local-kernels} we have for any $\sigma >0$
\begin{equation}\label{kernel-local}
|\varphi_k(x, y)|, |\psi_k(x, y)|, |\tilde\psi_k(x, y)|
\le c_\sigma |B(x, 2^{-k})|^{-1}\big(1+2^k\rho(x, y)\big)^{-\sigma}.
\end{equation}

The following lemma will be instrumental in this proof:


\begin{lemma}\label{lem:estim-psi-phi}
Let $\phi\in\cS$, $k\ge 0$, and $\sigma>0$.
Then
\begin{equation}\label{est-psi-phi}
\Big|\int_M \psi_k(x, y)\phi(y)d\mu(y)\Big|
\le c_\sigma 2^{-k(\NN-d)}(1+\rho(x, x_0))^{-\sigma},
\quad \forall x\in M,
\end{equation}
where $c_\sigma=c\cP_n(\phi)$ with $n \ge \max\{\sigma+d, K/2\}$ and
$c>0$ is a constant independent of $k, \sigma, \phi$.
\end{lemma}

\noindent
{\bf Proof.}
Let $m :=\NN/2$ (recall that $\NN$ is even).
Setting $g(\lambda):= \lambda^{-2m}\psi(\lambda)$ we have
$
L^{-m}\psi(2^{-k}\sqrt{L}) = 2^{-2km} g(2^{-k}\sqrt{L})
$
and hence, for $\phi\in\cS$,
\begin{align*}
\psi_k\phi(x)= 2^{-2km} g(2^{-k}\sqrt{L})L^m\phi(x)
= 2^{-kK}\int_M g(2^{-k}\sqrt{L})(x, y)L^m\phi(y) d\mu(y).
\end{align*}
Using (\ref{derivatives}), it readily follows that $g\in \cS(\R)$ and $g$ is real-valued and even.
Then by virtue of  Theorem~\ref{thm:S-local-kernels} for any $k\ge 0$
\begin{align}\label{ker-g}
|g(2^{-k}\sqrt{L})(x, y)|
&\le c_\sigma|B(x, 2^{-k})|^{-1}(1+2^k\rho(x, y))^{-\sigma-d}\\
&\le c_\sigma 2^{kd}|B(x, 1)|^{-1}(1+\rho(x, y))^{-\sigma-d},\notag
\end{align}
where we used (\ref{doubling}).
On the other hand, as $\phi\in \cS$ we have on account of (\ref{norm-S})
$|L^m\phi(x)|\le c(1+\rho(x, x_0))^{-\sigma-d}$ with $c=\cP_n(\phi)$, $n \ge \max\{\sigma+d, m\}$.
From this and (\ref{ker-g}) we infer
\begin{align*}
|\psi_k\phi(x)|
&\le c2^{-k(K-d)}|B(x,1)|^{-1}\int_M (1+\rho(x, y))^{-\sigma-d}(1+\rho(y, x_0))^{-\sigma-d}d\mu(y)\\
&\le c2^{-k(K-d)}(1+\rho(x, x_0))^{-\sigma}.
\end{align*}
Here for the last inequality we used (\ref{tech-2}).
$\qed$

\smallskip

In the following we will utilize the following assertion involving the grand maximal operator
$\cM_N$, defined in (\ref{def-grand-max3}):
Let $\phi\in \cS(\R)$ be admissible and assume $\cN_N(\phi)\le c$.
Then for any $f\in\cS'$, $k\in\bZ$, and $x\in M$
\begin{equation}\label{argument}
|\phi(2^{-k}\sqrt L)f(y)| \le c\cM_N(f)(x)
\quad\hbox{for all}\;\; y\in M \;\;{with}\;\; \rho(x, y)\le 2\tau 2^{-k},
\end{equation}
where $\tau >1$ is the constant from (\ref{supp-psi}).
This claim follows readily from (\ref{est-grand-max}).

\smallskip

%
Given $f\in H^p$, $0<p\le 1$, $f\ne 0$, we define
\begin{equation}\label{def-Omega-r}
\Omega_r:= \{x\in M: \cM_N(f)(x) >2^r\}, \quad r\in \bZ.
\end{equation}
Clearly,
$\Omega_{r+1} \subset \Omega_r$ and $M = \cup_{r\in \bZ} \Omega_r$.
The latter identity follows by $\cM_N(f)(x) >0$ $\forall x\in M$ due to $f\ne 0$.
Also, $\Omega_r$ is open since $\cM_N(f)(x)$ is lower semi-continuous.

It is easy to see that 
\begin{equation}\label{ineq-Omega}
\sum_{r\in \bZ} 2^{pr}|\Omega_r| \le c \int_M \cM_N(f)(x)^p d\mu(x) \le c\|f\|_{H^p}^p.
\end{equation}


We next show that
\begin{equation}\label{varphi-conv}
\|\varphi_j^2f\|_\infty \to 0 \quad\hbox{as}\quad j\to -\infty.
\end{equation}
Indeed, observe first that (\ref{ineq-Omega}) implies
$|\Omega_r|\le c2^{-pr}\|f\|_{H^p}^p$ for $r\in\bZ$.
%
Fix $r\in \bZ$. Then in light of (\ref{ball-low-est}) there exists $J>0$ s.t.
$\Omega_r \subset \{x\in M:\dist(x, \Omega_r^c) \le 2\cc 2^{J}\}$.
Hence for every $x\in M$ there exists $y\in\Omega_r^c$
such that $\rho(x, y) \le 2\cc 2^{J}$,
and by (\ref{argument}) we get
$$
|\varphi_j^2f(x)| \le c\cM_N(f)(y) \le c2^r
\quad\hbox{for $\;j \le -J$.}
$$
Hence
$\|\varphi_j^2f\|_\infty \le c2^r$ for $j \le -J$,
which implies (\ref{varphi-conv}).
From (\ref{main-repr}) and (\ref{varphi-conv}) it follows that
\begin{equation}\label{convergence}
f=\lim_{K\to\infty} \sum_{k=-\infty}^K \psi_k\tilde\psi_k f
\quad \hbox{in} \;\;\cS'.
\end{equation}

Assume $\Omega_r\ne \emptyset$ and write
\begin{equation}\label{def-E-rk}
E_{rk}:=\big\{x\in \Omega_r: \dist(x, \Omega_r^c) > 2\cc 2^{-k}\big\}
\setminus\big\{x\in \Omega_{r+1}: \dist(x, \Omega_{r+1}^c) > 2\cc 2^{-k}\big\}.
\end{equation}
%
By (\ref{ineq-Omega}) it follows that for any $r\in\bZ$ we have $|\Omega_r| <\infty$
and hence using (\ref{ball-low-est}) there exists $\kr\in \bZ$
such that $E_{r\kr}\ne\emptyset$ and $E_{rk}=\emptyset$ for $k<\kr$.
Note that $s_r \le s_{r+1}$.
We define
\begin{equation}\label{def-F-rr}
F_{r}(x):= \sum_{k \ge \kr} \int_{E_{rk}} \psi_k(x, y)\tilde\psi_kf(y) d\mu(y),
\quad x\in M, \; r\in \bZ,
\end{equation}
and in general
\begin{equation}\label{def-F-rk}
F_{r, \kl, \kll}(x):= \sum_{k = \kl}^{\kll} \int_{E_{rk}} \psi_k(x, y)\tilde\psi_kf(y) d\mu(y),
\quad \kr\le \kl \le \kll\le \infty.
\end{equation}
As will be shown in Lemma~\ref{lem:F-rr} below, the functions $F_r$ and $F_{r, \kl, \kll}$
are well defined and $F_r, F_{r, \kl, \kll}\in L^\infty$.

Observe that the fact that $\supp \psi(2^{-k}\sqrt L)(x, \cdot) \subset B(x, \cc 2^{-k})$
leads to the following conclusions:

(i) If $B(x, \cc 2^{-k}) \subset E_{rk}$ for some $x\in E_{rk}$, then
\begin{align}\label{local}
\int_{E_{rk}} \psi_k(x, y)\tilde\psi_kf d\mu(y)
= \int_{B(x, \cc 2^{-k})} \psi_k(x, y)\tilde\psi_k f(y) d\mu(y)
= \int_M \psi_k(x, y)\tilde\psi_kf(y) d\mu(y).
\end{align}

(ii) We have
\begin{equation}\label{supp-block}
\supp \Big(\int_{E_{rk}} \psi_k(\cdot, y)\tilde\psi_kf(y) d\mu(y)\Big)
\subset \big\{x: \dist(x,E_{rk})\le\cc 2^{-k}\big\}.
\end{equation}
On the other hand, clearly
$B(y,2\cc 2^{-k}) \cap \big(\Omega_r\setminus\Omega_{r+1}\big) \ne \emptyset$
for each $y\in E_{rk}$,
and $\cN_N(\tilde\psi) \le c$.
Therefore, see (\ref{argument}),
$|\tilde\psi_kf(y)|\le c2^r$ for $y\in E_{rk}$,
and using (\ref{kernel-local}) with $\sigma>d$ and (\ref{tech-1}) we get
\begin{equation}\label{norm-block}
\Big\|\int_{E} \psi_k(\cdot, y)\tilde\psi_kf(y) d\mu(y)\Big\|_\infty \le c2^r,
\quad \forall E \subset E_{rk}.
\end{equation}
Similarly,
\begin{equation}\label{norm-block-2}
\Big\|\int_{E} \varphi_k(\cdot, y)\varphi_k f(y) d\mu(y)\Big\|_\infty \le c2^r,
\quad \forall E \subset E_{rk}.
\end{equation}

We record some of the main properties of $F_r$ and $F_{r, \kl, \kll}$ in the following


\begin{lemma}\label{lem:F-rr}
$(a)$
We have
\begin{equation}\label{Erk-cover}
E_{rk}\cap E_{r'k}=\emptyset
\quad\hbox{if $r\ne r'\;\;$ and}
\quad M=\cup_{r\in \bZ} E_{rk}, \quad \forall k\in\bZ.
\end{equation}

$(b)$
There exists a constant $c>0$ such that for any $r\in\bZ$ and $s_r \le \kl \le \kll \le \infty$
\begin{equation}\label{est-F-rr}
\|F_r\|_\infty \le c2^r \quad\hbox{and}\quad \|F_{r, \kl, \kll}\|_\infty \le c2^r.
\end{equation}

$(c)$
The series in $(\ref{def-F-rr})$ and $(\ref{def-F-rk})$ $($if $ \kll=\infty)$ converge point-wise
and in distributional sense $($in $\cS'$$)$.

$(d)$
Also,
\begin{equation}\label{F-rr-zero}
F_r(x)=0\quad\hbox{for}\;\; x\in M\setminus\Omega_r, \;\; \forall r\in\bZ.
\end{equation}
\end{lemma}

\noindent
{\bf Proof.}
Identities (\ref{Erk-cover}) are obvious and
(\ref{F-rr-zero}) follows readily from the definition of $F_r$ and (\ref{supp-block}).

We next focus on the proof of the left-hand side estimate in (\ref{est-F-rr}).
The proof of the right-hand side estimate in (\ref{est-F-rr}) is the same; we omit it.
Assume $\Omega_{r+1} \ne \emptyset$;
the case when $\Omega_{r+1} = \emptyset$ is easier and will be omitted.
Set
$$ 
U_k=\big\{x\in \Omega_r: \dist(x, \Omega_r^c) > 2\cc 2^{-k}\big\},
\quad
V_k=\big\{x\in \Omega_{r+1}: \dist(x, \Omega_{r+1}^c) > 2\cc 2^{-k}\big\}.
$$ 
Then $E_{rk} = U_k\setminus V_k$, see (\ref{def-E-rk}).

From (\ref{supp-block}) we get $|F_r(x)|=0$ for
$x \in M\setminus \bigcup_{k\ge \kr}\{y: \dist(y, E_{rk})<\cc2^{-k}\}$.

We next estimate $|F_r(x)|$ for
$x \in \bigcup_{k\ge \kr} \{y: \dist(y, E_{rk})<\cc2^{-k}\}$.
Two cases present themselves here.

{\bf Case 1:}
$x\in \big[\cup_{k\ge \kr}\big\{y: \dist(y, E_{rk})<\cc2^{-k}\big\}\big]\cap\Omega_{r+1}$.
Then there exist $\nu, \ell\in\bZ$ such that
\begin{equation}\label{location-x}
x\in (U_{\ell+1}\setminus U_\ell)\cap (V_{\nu+1}\setminus V_\nu).
\end{equation}
Due to $\Omega_{r+1}\subset\Omega_r$ we have $V_k\subset U_k$, implying
$(U_{\ell+1}\setminus U_\ell)\cap (V_{\nu+1}\setminus V_\nu)=\emptyset$ if $\nu<\ell$.
We consider two subcases depending on whether $\nu\ge \ell+3$ or $\ell \le \nu \le \ell+2$.

(a) Let $\nu\ge \ell+3$.
We claim that (\ref{def-E-rk}) and (\ref{location-x}) yield
\begin{equation}\label{intersect-1}
B(x, \cc2^{-k}) \cap E_{rk} =\emptyset \quad\hbox{for}\quad k\ge \nu+2 \;\;\hbox{or}\;\; k\le \ell-1.
\end{equation}
Indeed, if $k\ge \nu+2$, then $E_{rk} \subset \Omega_r\setminus V_{\nu+2}$, which implies (\ref{intersect-1}),
while if $k\le \ell-1$, then $E_{rk} \subset U_{\ell-1}$, again implying (\ref{intersect-1}).

We also claim that
\begin{equation}\label{claim-1}
B(x, \cc2^{-k}) \subset E_{rk}\quad\hbox{for}\quad \ell+2\le k\le \nu-1.
\end{equation}
Indeed, if $\ell+2\le k\le \nu-1$, then
\begin{align*}
(U_{\ell+1}\setminus U_\ell)\cap (V_{\nu+1}\setminus V_\nu)
\subset (U_{k-1}\setminus U_\ell)\cap (V_{\nu+1}\setminus V_{k+1})
\subset U_{k-1}\setminus V_{k+1},
\end{align*}
which implies (\ref{claim-1}).

From (\ref{local})-(\ref{supp-block}) and (\ref{intersect-1})-(\ref{claim-1}) it follows that
\begin{align*}
F_{r}(x)&= \sum_{k =\ell}^{\nu+1} \int_{E_{rk}} \psi_k(x, y)\tilde\psi_kf(y) d\mu(y)
= \sum_{k =\ell}^{\ell+1} \int_{E_{rk}} \psi_k(x, y)\tilde\psi_kf(y) d\mu(y)\\
&+ \sum_{k =\ell+2}^{\nu-2} \int_M \psi_k(x, y)\tilde\psi_kf(y) d\mu(y)
+ \sum_{k =\nu-1}^{\nu+1} \int_{E_{rk}} \psi_k(x, y)\tilde\psi_kf(y) d\mu(y).
\end{align*}
However,
\begin{align*}
&\sum_{k =\ell+2}^{\nu-2} \int_M \psi_k(x, y)\tilde\psi_kf(y) d\mu(y)
= \sum_{k =\ell+2}^{\nu-2}\psi_k\tilde\psi(2^{-k}\sqrt L)f(x)\\
&\qquad\qquad
= \sum_{k =\ell+2}^{\nu-2} \big[\varphi_{k+1}^2f(x)-\varphi_k^2f(x)\big]
= \varphi_{\nu-1}^2f(x)-\varphi_{\ell+2}^2f(x)\\
&\qquad\qquad
=\int_{E_{r, \nu-1}} \varphi_{\nu-1}(x, y)\varphi_{\nu-1} f(y) dy
- \int_{E_{r, \ell+2}} \varphi_{\ell+2}(x, y)\varphi_{\ell+2}f(y) dy.
\end{align*}
Thus $F_r(x)$ is represented as the sum of at most seven integrals.
We estimate each of them using (\ref{norm-block})-(\ref{norm-block-2}) to obtain
$|F_r(x)| \le c2^r$.

(b) Let $\ell \le \nu \le \ell+2$.
We have
\begin{align*}
F_{r}(x)&= \sum_{k =\ell}^{\nu+1} \int_{E_{rk}} \psi_k(x, y)\tilde\psi_kf(y) d\mu(y)
= \sum_{k =\ell}^{\ell+3} \int_{E_{rk}} \psi_k(x, y)\tilde\psi_kf(y) d\mu(y).
\end{align*}
We use (\ref{norm-block}) to estimate each of these four integrals and obtain
again $|F_r(x)| \le c2^r$.

\smallskip


{\bf Case 2:} $x\in \Omega_r\setminus\Omega_{r+1}$.
Then there exists $\ell \in \bZ$ such that
$$
x\in (U_{\ell+1}\setminus U_\ell)\cap(\Omega_r\setminus\Omega_{r+1}).
$$
Just as in the proof of (\ref{intersect-1}) we have
$B(x, \cc2^{-k}) \cap E_{rk} =\emptyset$ for $k\le \ell-1$,
and as in the proof of (\ref{claim-1}) we have
$$ 
(U_{\ell+1}\setminus U_\ell)\cap (\Omega_r\setminus \Omega_{r+1})
\subset U_{k-1}\setminus V_{k+1}
\quad\hbox{for}\;\; k\ge \ell+2,
$$ 
which implies
$B(x, \cc2^{-k}) \subset E_{rk}$ for $k\ge \ell+2$.
We use these and (\ref{local})-(\ref{supp-block}) to obtain
\begin{align*}
F_{r}(x)
&= \sum_{k =\ell}^\infty \int_{E_{rk}} \psi_k(x, y)\tilde\psi_kf(y) d\mu(y)\\
&= \sum_{k =\ell}^{\ell+1} \int_{E_{rk}} \psi_k(x, y)\tilde\psi_kf(y) d\mu(y)
+ \sum_{k =\ell+2}^\infty \int_M \psi_k(x, y)\tilde\psi_kf(y) d\mu(y).
\end{align*}
For the last sum we have
\begin{align*}
&\sum_{k =\ell+2}^\infty \int_M \psi_k(x, y)\tilde\psi_kf(y) d\mu(y)\\
&\quad
= \lim_{\nu\to \infty}\sum_{k =\ell+2}^{\nu}\psi_k\tilde\psi_kf(x)
= \lim_{\nu\to \infty}\big[\varphi_{\nu+1}^2f(x) - \varphi_{\ell+2}^2f(x)\big]\\
&\quad
= \lim_{\nu\to \infty}\Big(\int_{E_{r, \nu+1}}
\varphi_{\nu+1}(x, y)\varphi_{\nu+1}f(y) d\mu(y)
- \int_{E_{r, \ell+2}} \varphi_{\ell+2}(x, y)\varphi_{\ell+2}f(y) d\mu(y)\Big).
\end{align*}
From the above and (\ref{norm-block})-(\ref{norm-block-2}) we enfer
$|F_r(x)| \le c2^r$.

\smallskip

The point-wise convergence in (\ref{def-F-rr}) follows from above
and we similarly establish the point-wise convergence in (\ref{def-F-rk}).

The convergence in distributional sense in (\ref{def-F-rr}) relies on the following assertion:
For every $\phi\in\cS$
\begin{equation}\label{conv-distr}
\sum_{k \ge \kr} |\langle g_{rk}, \phi\rangle| <\infty,
\quad\hbox{where}\quad g_{rk}(x):=\int_{E_{rk}} \psi_k(x, y)\tilde\psi_kf(y) d\mu(y).
\end{equation}
Here $\langle g_{rk}, \phi\rangle:= \int_M g_{rk} \overline{\phi} d\mu$.
To prove the above we need this estimate:
\begin{equation}\label{est-psi-f}
|\tilde\psi_k f(x)| 
\le\left\{
\begin{array}{ll}
c2^{kd/p}|B(x,1)|^{-1/p}\|f\|_{H^p}, & k\ge 0,\\
c2^{\eps k/p}|B(x,1)|^{-1/p}\|f\|_{H^p}, & k<0.
\end{array}
\right.
\end{equation}
Indeed, using (\ref{est-grand-max}) we get
\begin{align*}
|\tilde\psi_k f(x)|^p
&\le \inf_{y: \rho(x, y)\le 2^{-k}} \sup_{z: \rho(y, z)\le 2^{-k}}|\tilde\psi_k f(z)|^p
\le \inf_{y: \rho(x, y)\le 2^{-k}} c\cM_N(f)(y)^p\\
& \le c|B(x, 2^{-k})|^{-1}\int_{B(x, 2^{-k})}\cM_N(f)(y)^p d\mu(y)
\le c|B(x, 2^{-k})|^{-1}\|f\|_{H^p}^p.
\end{align*}
Then (\ref{est-psi-f}) follows by (\ref{doubling}) and (\ref{reverse-doubling-2}).

We now estimate $|\langle g_{rk}, \phi\rangle|$.
From (\ref{kernel-local}), (\ref{est-psi-f}), and the fact that $\phi\in\cS$ it readily follows that
$$
\int_{E_{rk}}\int_M |\psi_k(x, y)||\phi(x)||\tilde\psi_kf(y)| d\mu(y) d\mu(x) <\infty,
\quad k \ge s_r.
$$
Therefore, we can use Fubini's theorem, Lemma~\ref{lem:estim-psi-phi} (with $\sigma >2d$), and (\ref{est-psi-f})
to obtain for $k\ge 0$
\begin{align}\label{gk-phi}
|\langle g_{rk}, \phi\rangle|
& \le \int_{E_{rk}} \Big|\int_M \psi_k(x, y)\phi(x) d\mu(x)\Big||\tilde\psi_kf(y)| d\mu(y)\notag\\
& = \int_{E_{rk}} \Big|\int_M \psi_k(y, x)\phi(x) d\mu(x)\Big||\tilde\psi_kf(y)| d\mu(y)\\
& \le c2^{-k(\NN-d-d/p)}\|f\|_{H^p} \int_{E_{rk}}|B(y, 1)|^{-1/p}(1+\rho(y, x_0))^{-\sigma} d\mu(y).\notag
\end{align}
Here we also used that $\psi_k(y, x)=\psi_k(x, y)$ because $\psi$ is real-valued.
Further, from the non-collapsing condition (\ref{non-collapsing}) we have $|B(y, 1)|\ge c_2>0$
and using (\ref{tech-1}) we arrive at
$$
|\langle g_{rk}, \phi\rangle| \le c2^{-k(\NN-d-d/p)}\|f\|_{H^p} \quad\hbox{for $k\ge 0$.}
$$
This implies (\ref{conv-distr}) due to $\NN > d+d/p$.

Write $G_\ell := \sum_{k=s_r}^\ell g_{rk}$.
From the above proof of (b) and (\ref{est-F-rr}) we infer that $G_\ell(x) \to F_r(x)$ as $\ell \to \infty$
for $x\in M$ and $\|G_\ell\|_\infty \le c2^r< \infty$ for $\ell \ge s_r$.
On the other hand, from (\ref{conv-distr}) it follows that the series $\sum_{k\ge s_r} g_{rk}$
converges in distributional sense.
By applying the dominated convergence theorem one easily concludes that $F_r=\sum_{k\ge s_r} g_{rk}$
with the convergence in distributional sense.
$\qed$

\smallskip

For convenience, we define $F_r:=0$ whenever $\Omega_r=\emptyset$, $r\in \bZ$.

Observe that from (\ref{Erk-cover}) it follows that
\begin{equation}\label{rep-psi-psi}
\psi_k\tilde\psi_kf(x)
= \int_M \psi_k(x, y)\tilde\psi_kf(y)d\mu(y)
= \sum_{r\in\bZ}\int_{E_{rk}} \psi_k(x, y)\tilde\psi_kf(y)d\mu(y)
\end{equation}
and using $(\ref{convergence})$ and the definition of $F_r$ in $(\ref{def-F-rr})$ we arrive at
\begin{equation}\label{represent}
f = \sum_{r\in\bZ} F_r \;\; \hbox{in}\;\; \cS',
\;\;\hbox{that is,}\quad
\langle f, \phi\rangle = \sum_{r\in\bZ} \langle F_r, \phi\rangle,
\quad\forall \phi\in\cS,
\end{equation}
where the last series converges absolutely.
We next give the needed justification of identity (\ref{represent}).

From (\ref{convergence}), (\ref{def-F-rr}), (\ref{rep-psi-psi}),
and the notation from (\ref{conv-distr}) we obtain for $\phi\in\cS$
\begin{align*}
\langle f, \phi\rangle
= \sum_{k}\langle\psi_k\tilde\psi_kf, \phi \rangle
= \sum_{k}\sum_{r}\langle g_{rk}, \phi\rangle
= \sum_{r}\sum_{k}\langle g_{rk}, \phi\rangle
= \sum_{r}\langle F_r, \phi\rangle.
\end{align*}
Clearly, to justify the above identities it suffices to show that
$
\sum_{k}\sum_{r} |\langle g_{rk}, \phi\rangle| <\infty.
$
We split this sum into two:
$
\sum_{k}\sum_{r}\cdots
= \sum_{k\ge 0}\sum_{r} \cdots + \sum_{k<0}\sum_{r} \cdots
=: \Sigma _1+ \Sigma_2.
$
To estimate $\Sigma_1$ we use (\ref{gk-phi}) and that $|B(y, 1)|\ge c_2$, by (\ref{non-collapsing}).
We obtain
\begin{align*}
\Sigma_1
&\le c\|f\|_{H^p}  \sum_{k\ge 0}2^{-k(\NN-d-d/p)}\sum_{r}
\int_{E_{rk}}|B(y, 1)|^{-1}(1+\rho(y, x_0))^{-\sigma} d\mu(y)\\
&\le c\|f\|_{H^p}  \sum_{k\ge 0}2^{-k(\NN-d-d/p)}
\int_M|B(y, 1)|^{-1}(1+\rho(y, x_0))^{-\sigma} d\mu(y)
\le c\|f\|_{H^p}.
\end{align*}
Here we also used that $K> d+d/p$, $\sigma>2d$, and (\ref{tech-1}).

We estimate $\Sigma_2$ in a similar manner, using (\ref{est-psi-f}),
again (\ref{non-collapsing}),
and the fact that $\int_M|\psi_k(x, y)|d\mu(y) \le c<\infty$
and $|\phi(x)| \le c(1+\rho(x, x_0))^{-d-1}$.
We get
\begin{align*}
\Sigma_2
&\le c\|f\|_{H^p}  \sum_{k<0}2^{\eps k/p}\sum_{r}
\int_{E_{rk}}\int_M |\psi_k(x, y)|d\mu(y)|B(x, 1)|^{-1/p}|\phi(x)| d\mu(x)\\
&\le c\|f\|_{H^p}  \sum_{k<0}2^{\eps k/p}\int_M |B(x, 1)|^{-1}(1+\rho(x, x_0))^{-d-1} d\mu(x)
\le c\|f\|_{H^p}.
\end{align*}
The above estimates of $\Sigma_1$ and $\Sigma_2$ imply
$
\sum_{k}\sum_{r} |\langle g_{rk}, \phi\rangle| <\infty,
$
which completes the justification of (\ref{represent}).


\medskip

We next break each function $F_r$ into atoms.
To this end we need a Whitney type cover for $\Omega_r$.


\begin{lemma}\label{lem:Whitney}
Assume $\Omega$ is an open proper subset of $M$ and
let $\rho(x):= \dist (x, \Omega^c)$.
Then there exist a constant $K >0$ $($$K = 70^d c_0^2$ will do$)$
and a sequence of points $\{\xi_j\}_{j\in\bN}$ in $\Omega$
with the following properties,
where $\rho_j:= \dist (\xi_j, \Omega^c)$:

\medskip

$(a)$ $\Omega = \cup_{j\in \bN} B(\xi_j, \rho_j/2)$.

\medskip

$(b)$ $\{B(\xi_j, \rho_j/5)\}$ are disjoint.

%
%
%

\medskip

$(c)$ If $B\big(\xi_j, \frac{3\rho_j}{4}\big)\cap B\big(\xi_\nu, \frac{3\rho_\nu}{4}\big)\ne \emptyset$,
then  $7^{-1}\rho_\nu\le \rho_j \le 7\rho_\nu$.

\medskip

$(d)$ For every $j\in\bN$ there are at most $K$ balls $B\big(\xi_\nu, \frac{3\rho_\nu}{4}\big)$
intersecting $B\big(\xi_j, \frac{3\rho_j}{4}\big)$.

\end{lemma}

Variants of this simple lemma are well known and frequently used.
To prove it one simply selects $\{B(\xi_j, \rho(\xi_j)/5)\}_{j\in\bN}$
to be a maximal disjoint subcollection of $\{B(x, \rho(x)/5)\}_{x\in\Omega}$
and then properties (a)-(d) follow readily, see \cite{Stein}, pp. 15-16.
For completeness we give its proof in the appendix.


\smallskip

We apply Lemma~\ref{lem:Whitney} to each set $\Omega_r\ne \emptyset$. 
Fix $r\in \bZ$ and assume $\Omega_r\ne \emptyset$.
Denote by $B_j:= B(\xi_j, \rho_j/2)$, $j=1, 2, \dots$,
the balls given by Lemma~\ref{lem:Whitney}, applied to $\Omega_r$,
with the additional assumption that these balls are ordered so that
$\rho_1 \ge \rho_2 \ge \cdots$.
We will adhere to the notation from Lemma~\ref{lem:Whitney}.
We will also use the more compact notation
$\cB_r:=\{B_j\}_{j\in\bN}$ for the set of balls covering $\Omega_r$.

For each ball $B\in \cB_r$ and $k\ge \kr$ we define
\begin{equation}\label{def-E-B}
E_{rk}^B:=E_{rk} \cap \big\{x: \dist(x, B) < 2\cc 2^{-k}\big\}
\quad\hbox{if}\quad B\cap E_{rk}\ne \emptyset
\end{equation}
and set
$E_{rk}^B:=\emptyset$ if $B\cap E_{rk}= \emptyset$.

We also define, for $\ell=1, 2, \dots$,
\begin{equation}\label{def-R-B}
R_{rk}^{B_\ell}:= E_{rk}^{B_\ell}\setminus \cup_{\nu>\ell}E_{rk}^{B_\nu},
\end{equation}
\begin{equation}\label{def-F-B}
F_{B_\ell}(x):= \sum_{k \ge \kr} \int_{R_{rk}^{B_\ell}}
\psi_k(x, y)\tilde\psi_kf(y) d\mu(y),
\quad\hbox{and} \quad G_{B_\ell}:=L^{-n}F_{B_\ell}.
\end{equation}


\begin{lemma}\label{lem:FrkB}
We have for any $k\ge \kr$
\begin{equation}\label{FrkB1}
E_{rk}=\cup_{\ell\ge 1} R_{rk}^{B_\ell}
\quad\hbox{and}\quad R_{rk}^{B_\ell}\cap R_{rk}^{B_m} = \emptyset
\quad\hbox{if} \;\;\ell\ne m.
\end{equation}
Hence
\begin{equation}\label{decomp-Fr}
F_r= \sum_{B\in\cB_r} F_{B}\quad\hbox{$($convergence in $\cS'$$)$.}
\end{equation}
Furthermore, the series in $(\ref{def-F-B})$ converges point-wise and in $\cS'$, and
there exists a constant $\cs >0$ such that for every $\ell \ge 1$
\begin{equation}\label{supp-FB}
\supp F_{B_\ell} \subset 7B_\ell,
\end{equation}
\begin{equation}\label{est-FB}
\|F_{B_\ell}\|_\infty \le \cs 2^r
\quad\hbox{and}\quad
\|L^m G_{B_\ell}\|_\infty \le \cs 2^r \rho_\ell^{2(n-m)} \quad\hbox{for}\;\; m= 0, 1, \dots, n.
\end{equation}
\end{lemma}

\noindent
{\bf Proof.}
By Lemma~\ref{lem:Whitney} we have $\Omega_r=\cup_{\ell\in\bN} B_\ell$
and then (\ref{FrkB1}) is immediate from (\ref{def-E-B}) and (\ref{def-R-B}).

Fix $\ell \ge 1$. Observe that using Lemma~\ref{lem:Whitney} we have
$B_\ell \subset \{x: \dist(x, \Omega_r^c) <  2\rho_\ell\}$
and hence $E_{rk}^{B_\ell}:=\emptyset$ if $2\cc 2^{-k} \ge 2\rho_\ell$.
Define $\kkl: = \min \{k: \cc 2^{-k} < \rho_\ell\}$.
Hence $\rho_\ell/2 \le \cc 2^{-\kkl} <\rho_\ell$.
Consequently,
\begin{equation}\label{def-F-B-2}
F_{B_\ell}(x)= \sum_{k \ge \kkl} \int_{R_{rk}^{B_\ell}}
\psi_k(x, y)\tilde\psi_kf(y) d\mu(y).
\end{equation}
Clearly
$\supp F_{B_\ell} \subset B\big(\xi_\ell, (7/2)\rho_\ell\big)=7B_\ell$,
which confirms (\ref{supp-FB}).

To prove the left-hand side estimate in (\ref{est-FB}) we need some preparation.


\begin{lemma}\label{lem:FS}
For an arbitrary set $S\subset M$ let
$
S_k:= \big\{x\in M: \dist(x, S) < 2\cc 2^{-k}\big\}
$
and set
\begin{equation}\label{def-FS}
F_S(x):= \sum_{k \ge \kl} \int_{E_{rk}\cap S_k} \psi_k(x, y)\tilde\psi_kf(y) d\mu(y)
\end{equation}
for some $\kl \ge \kr$.
Then
$
\|F_S\|_\infty \le c2^r,
$
where $c>0$ is a constant independent of $S$ and $\kl$.
Moreover, the above series converges in $\cS'$.
\end{lemma}

\noindent
{\bf Proof.}
From (\ref{supp-block}) it follows that $F_S(x)=0$ if $\dist (x, S) \ge 3\cc 2^{-\kl}$.

Let $x\in S$. Clearly,
$B(x, \cc2^{-k})\subset S_k$ for every $k$ and hence
\begin{align*}
F_S(x)
&=\sum_{k \ge \kl} \int_{E_{rk}\cap B(x, \cc2^{-k})} \psi_k(x, y)\tilde\psi_kf(y) d\mu(y)\\
&=\sum_{k \ge \kl} \int_{E_{rk}} \psi_k(x, y)\tilde\psi_kf(y) d\mu(y)= F_{r,\kl}(x).
\end{align*}
On account of Lemma~\ref{lem:F-rr} (b) we obtain
$|F_S(x)|= |F_{r,\kl}(x)|\le c2^r$.

Consider the case when $x\in S_\ell\setminus S_{\ell+1}$ for some $\ell \ge \kl$.
Then $B(x, \cc2^{-k}) \subset S_k$ if $\kl\le k\le \ell-1$
and $B(x, \cc2^{-k}) \cap S_k=\emptyset$ if $k\ge \ell+2$.
Therefore,
\begin{align*}
F_S(x) &=\sum_{k = \kl}^{\ell-1} \int_{E_{rk}} \psi_k(x, y)\tilde\psi_kf(y) d\mu(y)
+ \sum_{k = \ell}^{\ell+1} \int_{E_{rk}\cap S_k} \psi_k(x, y)\tilde\psi_kf(y) d\mu(y)\\
&= F_{r, \kl, \ell-1}(x)+ \sum_{k = \ell}^{\ell+1} \int_{E_{rk}\cap S_k} \psi_k(x, y)\tilde\psi_kf(y) d\mu(y),
\end{align*}
where we used the notation from (\ref{def-F-rk}).
By Lemma~\ref{lem:F-rr} (b) and (\ref{norm-block}) 
it follows that $|F_S(x)|\le c2^r$.

We finally consider the case when $2\cc 2^{-\kl} \le \dist (x, S) < 3\cc 2^{-\kl}$.
Then we have
$
F_S(x) = \int_{E_{r\kl}\cap S_\kl} \psi_\kl(x, y)\tilde\psi_\kl f(y) d\mu(y)
$
and the estimate $|F_S(x)|\le c2^r$ is immediate from (\ref{norm-block}).

The convergence in $\cS'$ in (\ref{def-FS}) is established just as in the proof of Lemma~\ref{lem:F-rr}.
$\qed$

\medskip
With $\ell\ge 1$ being fixed, we
let $\{B_j: j\in \cJ\}$ denote the set of all balls $B_j= B(\xi_j, \rho_j/2)$
such that $j>\ell$ and
$$
B\Big(\xi_j, \frac{3\rho_j}{4}\Big) \cap B\Big(\xi_\ell, \frac{3\rho_\ell}{4}\Big) \ne \emptyset.
$$
By Lemma~\ref{lem:Whitney} it follows that $\# \cJ \le K$
and $7^{-1} \rho_\ell \le \rho_j \le \rho_\ell$ for $j\in\cJ$.
Define
\begin{equation}\label{def-k1}
k_1:= \min\Big\{k: 2\cc2^{-k} < 4^{-1}\min\big\{\rho_j: j\in\cJ \cup\{\ell\}\big\} \Big\}.
\end{equation}
From this definition and  $\cc2^{-\kkl} <\rho_\ell$ we infer
\begin{equation}\label{k-kkl}
2\cc2^{-k_1} \ge 8^{-1} \min\big\{\rho_j: j\in\cJ \cup\{\ell\}\big\}
> 8^{-2}\rho_\ell > 8^{-2}\cc2^{-\kkl}
\; \Longrightarrow \; k_1 \le \kkl+7.
\end{equation}
Clearly, from (\ref{def-k1})
\begin{equation}\label{BB-1}
B(\xi_j, \rho_j/2+2\cc2^{-k}) \subset B\big(\xi_j, 3\rho_j/4\big),
\quad \forall k\ge k_1, \;\; \forall j\in \cJ \cup\{\ell\}.
\end{equation}
Denote
$S:= \cup_{j\in\cJ} B_j$ and
$\tilde{S}:= \cup_{j\in\cJ} B_j\cup B_\ell = S\cup B_\ell$.
As in Lemma~\ref{lem:FS} we set
$$
S_k:= \big\{x\in M: \dist(x, S) < 2\cc 2^{-k}\big\}
\quad\hbox{and}\quad
\tilde{S}_k:= \big\{x\in M: \dist(x, \tilde{S}) < 2\cc 2^{-k}\big\}.
$$
It readily follows from the definition of $k_1$ in (\ref{def-k1}) and (\ref{def-R-B}) that
\begin{equation}\label{R-EE}
R_{rk}^{B_\ell} := E_{rk}^{B_\ell}\setminus \cup_{\nu>\ell} E_{rk}^{B_\nu}
= \big(E_{rk}\cap \tilde{S}_k\big) \setminus \big(E_{rk}\cap S_k\big)
\quad\hbox{for}\quad k\ge k_1.
\end{equation}
Denote
\begin{align*}
F_S(x)&:= \sum_{k \ge k_1} \int_{E_{rk}\cap S_k} \psi_k(x, y)\tilde\psi_kf(y) d\mu(y),\quad\hbox{and}\\
F_{\tilde{S}}(x)&:= \sum_{k \ge k_1} \int_{E_{rk}\cap \tilde{S}_k} \psi_k(x, y)\tilde\psi_kf(y) d\mu(y).
\end{align*}
From (\ref{R-EE}) and the fact that $S \subset \tilde{S}$ it follows that
\begin{align*}
F_{B_\ell}(x) = F_{\tilde{S}}(x) - F_S(x)
+ \sum_{\kkl\le k < k_1} \int_{R_{rk}^{B_\ell}} \psi_k(x, y)\tilde\psi_kf(y) d\mu(y).
\end{align*}
By Lemma~\ref{lem:FS} we get
$\|F_S\|_\infty \le c2^r$ and $\|F_{\tilde{S}}\|_\infty \le c2^r$.
On the other hand from (\ref{k-kkl}) we have $k_1-\kkl \le 7$.
We estimate each of the (at most $7$) integrals above using (\ref{norm-block}) to conclude that
$\|F_{B_\ell}\|_\infty \le c2^r$.

The convergence in (\ref{def-F-B}) and (\ref{decomp-Fr}) is handled as in the proof of Lemma~\ref{lem:F-rr}.

\smallskip

It remains to prove that
$\|L^m G_{B_\ell}\|_\infty \le c2^r\rho_\ell^{2(n-m)}$, $0\le m <n$,
which is the second estimate in (\ref{est-FB}).
By definition
\begin{align*}
L^mG_{B_\ell}(x)
&:= L^{-(n-m)}F_{B_\ell}(x)\\
&= \sum_{k\ge \kkl} \int_{R_{rk}^{B_\ell}}
L^{-(n-m)}\psi(2^{-k}\sqrt L)(x, y)\tilde\psi(2^{-k}\sqrt L)f(y)d\mu(y).
\end{align*}
Let $g(\lambda):= \lambda^{-2(n-m)}\psi(\lambda)$.
Then $L^{-(n-m)}\psi(2^{-k}\sqrt L) = 2^{-2k(n-m)}g(2^{-k}\sqrt L)$.
From the definition of $\psi$, (\ref{derivatives}), and the fact that $K > 2n$ it follows that
$g\in \cS(\R)$ and $g$ is real-valued and even.
Therefore, by Theorem~\ref{thm:S-local-kernels} for any $k\ge 0$ and $\sigma>d$
$$
|L^{-(n-m)}\psi(2^{-k}\sqrt L)(x, y)|
\le \frac{c_\sigma 2^{-2k(n-m)}}{|B(x, 2^{-k})|(1+2^k\rho(x, y))^{\sigma}}.
$$
On the other hand by (\ref{norm-block})
$|\tilde\psi(2^{-k}\sqrt L)f(y)|\le c2^r$ for $y\in R_{rk}^{B_\ell}\subset E_{rk}$.
Putting the above together we get
\begin{align*}
&\Big|\int_{R_{rk}^{B_\ell}} L^{-(n-m)}\psi(2^{-k}\sqrt L)(x, y)\tilde\psi(2^{-k}\sqrt L)f(y)d\mu(y)\Big|\\
&\qquad\qquad\qquad
\le c2^r 2^{-2k(n-m)}\int_M \frac{d\mu(y)}{|B(x, 2^{-k})|(1+2^k\rho(x, y))^{\sigma}}
\le c2^r 2^{-2k(n-m)}.
\end{align*}
Hence,
$
\|L^m G_{B_\ell}\|_\infty \le c2^r\sum_{k\ge \kkl}2^{-2k(n-m)}
\le c2^r2^{-2\kkl(n-m)} \le c2^r \rho_\ell^{2(n-m)}
$
as claimed.
This completes the proof of Lemma~\ref{lem:FrkB}.
$\qed$

\smallskip

We are now in a position to complete the proof of Theorem~\ref{thm:hardy}.
For every ball $B\in \cB_r$, $r\in\bZ$, provided $\Omega_r\ne \emptyset$,
we define $B^\star:= 7B$,
$$
a_B(x):= \cs^{-1}|B^\star|^{-1/p}2^{-r}F_{B}(x),
\quad
b_B(x):= \cs^{-1}|B^\star|^{-1/p}2^{-r}G_{B}(x),
$$
and
$\lambda_B:= \cs|B^\star|^{1/p}2^r$,
where $\cs>0$ is the constant from (\ref{est-FB}).
By (\ref{supp-FB}) we have $\supp a_B\subset B^\star$,
and by (\ref{est-FB})
$$
\|a_B\|_\infty \le \cs^{-1}|B^\star|^{-1/p}2^{-r}\|F_{B}\|_\infty \le |B^\star|^{-1/p}.
$$
From (\ref{def-F-B}) it follows that $L^nb_B=a_B$
and assuming that $B=B(\xi_\ell, \rho_\ell/2)$ we obtain using (\ref{est-FB})
$$
\|L^mb_B\|_\infty
\le \cs^{-1}|B^\star|^{-1/p}2^{-r}\|L^mG_{B}\|_\infty
\le \rho_\ell^{2(n-m)}|B^\star|^{-1/p}
\le r_{B^\star}^{2(n-m)}|B^\star|^{-1/p}.
$$
Therefore, each $a_B$ is an atom for $H^p$.

We set $\cB_r:=\emptyset$ 
if $\Omega_r= \emptyset$.
Now, from above, (\ref{represent}), and Lemma~\ref{lem:FrkB} we infer
$$
f=\sum_{r\in\bZ} F_r
=\sum_{r\in\bZ}\sum_{B\in \cB_r}F_{B}
=\sum_{r\in\bZ}\sum_{B\in \cB_r}\lambda_B a_B,
$$
where the convergence is in $\cS'$, and
$$
\sum_{r\in\bZ}\sum_{B\in \cB_r}|\lambda_B|^p
\le c\sum_{r\in\bZ} 2^{pr}\sum_{B\in \cB_r}|B|
= c\sum_{r\in\bZ} 2^{pr}|\Omega_r|
\le c \|f\|_{H^p}^p,
$$
which is the claimed atomic decomposition of $f\in H^p$.
Above we used (\ref{ineq-Omega}) and that $|B^\star| = |7B| \le c_0 7^d |B|$.
$\qed$

\subsection{Proof of the embedding \boldmath $H^p_A \subset H^p$ in the noncompact case}
\label{subsec:embed-HpA-Hp}

We next show that if $f\in H^p_A$, then $f\in H^p$ and
$\|f\|_{H^p} \le c\|f\|_{H^p_A}$.
To this end we need the following


\begin{lemma}\label{lem:H-LP-atom}
For any atom $a$ and $0<p \le 1$, we have
\begin{equation}\label{F-atom}
\|a\|_{H^p} \le c<\infty.
\end{equation}
\end{lemma}

\noindent
{\bf Proof.}
Let $a(x)$ be an atom in the sense of Definition~\ref{def-atoms-hom}
and suppose $\supp a \subset B$, $B=B(z, r)$, and
$a=L^n b$ for some $b\in D(L^n)$, $\supp b\subset B$, and $\|b\|_\infty \le r^{2n} |B|^{-1/p}$.

Suppose $\varphi\in C^\infty_0(\R)$, $\varphi$ is real-valued and even,
$\supp \varphi \subset [-1,1]$,
$\varphi(0)=1$, and $\varphi^{(\nu)}(0)=0$ for $\nu \ge 1$.
By Theorem~\ref{thm:S-local-kernels},
applied with $f(\lambda)= \varphi(\lambda)$ and $f(\lambda)= \lambda^{2n}\varphi(\lambda)$,
it follows that $\varphi(t\sqrt{L})$ and $L^n\varphi(t\sqrt{L})$
are kernel operators with kernels satisfying the following inequalities for any $\sigma >0$
\begin{align}
|\varphi(t\sqrt{L})(x, y)| &\le c_\sigma |B(x, t)|^{-1}(1+ t^{-1}\rho(x, y))^{-\sigma} \label{kernel-1},\\
|L^n\varphi(t\sqrt{L})(x, y)| &\le c_\sigma t^{-2n}|B(x, t)|^{-1}(1+ t^{-1}\rho(x, y))^{-\sigma} \label{kernel-2}.
\end{align}
We choose $\sigma$ so that $\sigma > d/p+2d$.

We need estimate $|\varphi(t\sqrt{L})a(x)|$.
Observe first that using (\ref{tech-1}) we have
\begin{align}\label{atom-est-1}
|\varphi(t\sqrt{L})a(x)| \le \int_M \frac{|a(y)|}{|B(x, t)|(1+ t^{-1}\rho(x, y))^{\sigma}} d\mu(y)
\le c|B|^{-1/p},
\quad x\in 2B.
\end{align}

To estimate $|\varphi(t\sqrt{L})a(x)|$ for $x\in M\setminus 2B$ we consider two cases:

{\bf Case 1:} $0<t\le r$. Let $x\in M\setminus 2B$ and $y\in B$. From (\ref{doubling}) and (\ref{D2})
it readily follows that
$$
|B| 
\le c_0\Big(\frac{r}{t}\Big)^d|B(z, t)|
\le c_0^2\Big(\frac{r}{t}\Big)^d\Big(1+ \frac{\rho(x,z)}{t}\Big)^d  |B(x, t)|
\le c_0^2\Big(1+ \frac{\rho(x,z)}{t}\Big)^{2d}  |B(x, t)|,
$$
where we used that $\rho(x, z) \ge r$.
Combining this with (\ref{kernel-1}) and the obvious inequality
$\rho(x, z) \le \rho(x, y)+\rho(y, z) \le 2\rho(x, y)$
we obtain
$$
|\varphi(t\sqrt{L})(x, y)|
\le c_\sigma |B(x, t)|^{-1}(1+ t^{-1}\rho(x, y))^{-\sigma}
\le c |B|^{-1}(1+ t^{-1}\rho(x, z))^{-\sigma+2d}.
$$
In turn, this leads to
\begin{align*}
|\varphi(t\sqrt{L})a(x)|
&= \Big|\int_B \varphi(t\sqrt{L})(x, y)a(y)d\mu(y)\Big|\\
&\le \frac{c|B|^{-1-1/p}}{(1+ t^{-1}\rho(x, z))^{\sigma-2d}}
\int_B 1d\mu(y)
= \frac{c|B|^{-1/p}}{(1+ t^{-1}\rho(x, z))^{\sigma-2d}}.
\end{align*}
From this and (\ref{atom-est-1}) we infer
\begin{align}\label{norm-varphi-a}
\|\varphi(t\sqrt{L})a\|_{L^p}^p
&= \|\varphi(t\sqrt{L})a\|_{L^p(2B)}^p+ \|\varphi(t\sqrt{L})a\|_{L^p(M\setminus 2B)}^p\notag\\
&\le c\int_{2B}|B|^{-1} d\mu(x)
+ c\int_{M}\frac{|B|^{-1}d\mu(x)}{(1+ t^{-1}\rho(x, z))^{(\sigma-2d)p}}\\
& \le c'+c|B|^{-1} |B(z, t)| \le c.\notag
\end{align}
Here we used that $(\sigma-2d)p > d$ and (\ref{tech-1}).

{\bf Case 2:} $t > r$. Let $x\in M\setminus 2B$ and $y\in B$.
Using (\ref{D2}) we obtain
$$
|B| =|B(z, r)| \le |B(z, t)|
\le c_0\big(1+ \rho(x,z)/t\big)^{d}  |B(x, t)|
$$
and as before $\rho(x, z)\le 2\rho(x, y)$.
These coupled with (\ref{kernel-2}) lead to
$$
|L^n\varphi(t\sqrt{L})(x, y)| \le c (r/t)^{2n}|B|^{-1}(1+ t^{-1}\rho(x, z))^{-\sigma+d}.
$$
This and $\|b\|_\infty \le r^{2n} |B|^{-1/p}$ imply
\begin{align*}
|\varphi(t\sqrt{L})a(x)|
&= \Big|\int_B L^n\varphi(t\sqrt{L})(x, y)b(y)d\mu(y)\Big|\\
& \le \frac{c(r/t)^{2n}|B|^{-1-1/p}}{(1+ t^{-1}\rho(x, z))^{\sigma-2d}}
\int_B 1d\mu(y)
= \frac{c(r/t)^{2n}|B|^{-1/p}}{(1+ t^{-1}\rho(x, z))^{\sigma-2d}}.
\end{align*}
We use this and (\ref{atom-est-1}) to obtain
\begin{align*}
\|\varphi(t\sqrt{L})a\|_{L^p}^p
&= \|\varphi(t\sqrt{L})a\|_{L^p(2B)}^p+ \|\varphi(t\sqrt{L})a\|_{L^p(M\setminus 2B)}^p\\
&\le c\int_{2B}|B|^{-1} d\mu(x)
+ c\int_{M}\frac{(r/t)^{2np}|B|^{-1}d\mu(x)}{(1+ t^{-1}\rho(x, z))^{(\sigma-2d)p}}\\
& \le c'+c(r/t)^{2np}|B|^{-1} |B(z, t)|
\le c'+cc_0(r/t)^{2np}(t/r)^d\\
&= c'+cc_0(r/t)^{2np-d} \le c.
\end{align*}
Here we used that
$|B(z, t)| \le c_0(t/r)^d|B(z, r)|$
by (\ref{doubling}) and that $n\ge d/2p$.
In light of Theorem~\ref{thm:Hp} the above and (\ref{norm-varphi-a}) yield (\ref{F-atom}).
$\qed$

\smallskip

We are now prepared to complete the proof of the embedding $H^p_A \subset H^p$.
Assume that $f\in H^p_A$.
Then (see Definition~\ref{def-atoms-hom}) there exist atoms $\{a_k\}_{k\ge 1}$
and coefficients $\{\lambda_k\}_{k\ge 1}$ such that
$f=\sum_k{\lambda_k a_k}$ (convergence in $\cS'$) and $\sum_k{|\lambda_k|^p} \le 2\|f\|^p_{H^p_A}$.

Let
$\varphi\in C^\infty_0(\R)$ be real-valued, $\supp \varphi \subset [-1,1]$,
$\varphi(0)=1$, and $\varphi^{(\nu)}(0)=0$ for $\nu \ge 1$.
Then
$$
\varphi(t\sqrt L) f(x) = \sum_{k=1}^\infty \lambda_k \varphi(t\sqrt L)a_k(x), \quad x\in M, \;\;t>0,
$$
and hence
$$
\sup_{t>0}|\varphi(t\sqrt L) f(x)|
\le \sum_{k=1}^\infty |\lambda_k|\sup_{t>0} |\varphi(t\sqrt L)a_k(x)|,
$$
which is the same as
$
M(f;\varphi)(x) \le \sum_{k=1}^\infty |\lambda_k|M(a_k;\varphi)(x).
$
Therefore, for $ 0<p\le 1$
$$
\|M(f;\varphi)\|_{L^p}^p \le \sum_{k=1}^\infty |\lambda_k|^p \|M(a_k;\varphi)\|_{L^p}^p
\le c\sum_{k=1}^\infty |\lambda_k|^p \le c\|f\|_{H^p_A}^p.
$$
On account of Theorem~\ref{thm:Hp} this implies
$\|f\|_{H^p} \le c\|f\|_{H^p_A}$.
$\qed$

\subsection{Proof of Theorem~\ref{thm:hardy} in the compact case}
\label{subsec:Hp=HpA-compact}

We proceed quite similarly as in the noncompact case. 
Therefore, we will only indicate the modifications that need to be made.

To prove the embedding $H^p \subset H^p_A$
assume $f\in H^p$, $0<p\le 1$.
Let $\varphi\in \cS(\R)$ be just as in the proof in the noncompact case.
Instead of (\ref{convergence}) we use this representation of $f$ (see (\ref{main-repr})):
\begin{equation}\label{main-repr-inh}
f= \varphi_j^2f + \sum_{k=j+1}^\infty \psi_k\tilde\psi_k f
=: f_0+ f_1
\quad\hbox{(convergence in $\;\cS'$)},
\end{equation}
where $j$ is the maximal integer such that $B(x_0, 2^{-j})=M$,
and $\varphi_j$, $\psi_k$ and $\tilde\psi_k$ are as in (\ref{brief-notation}).
For the decomposition of $f_1$ we just repeat the proof from \S\ref{subsec:embed-Hp-HpA}.
On the other hand, as in (\ref{argument}) we have
$|\varphi_j^2f(x)| \le c\cM(f)(y)$, $\forall x, y\in M$,
and hence
$$
\|\varphi_j^2f\|_\infty \le c|M|^{-1/p}\|\cM(f)(y)\|_{L^p} \le c_*|M|^{-1/p}\|f\|_{H^p}.
$$
We define the outstanding atom $A$ (see (\ref{def-A})) by
$A(x):= c_*^{-1}\|f\|_{H^p}^{-1}\varphi_j^2f(x)$ and set $\lambda_A:= c_*\|f\|_{H^p}$.
Clearly, $\|A\|_\infty \le |B|^{-1/p}$
and
$\lambda_A A=\varphi_j^2f=f_0$.
Thus we arrive at the claimed atomic decomposition of $f$.

\smallskip

The proof of the embedding $H^p_A \subset H^p$
runs in the foot steps of the proof in the noncompact case from \S\ref{subsec:embed-HpA-Hp}.
We only have to show in addition the estimate
$\|A\|_{H^p} \le c <\infty$
for any outstanding atom $A$ as in (\ref{def-A}).
But, this estimate follows readily from estimate (\ref{atom-est-1}) applied to $A$.
$\qed$

\section{Decomposition of Hardy spaces via square functions}\label{sec:littlewood-paley}
\setcounter{equation}{0}

To put our study of Hardy spaces in prospective we bring here some relevant results.
In \cite{DKKP} we showed that in an inhomogeneous setting
the atomic Hardy spaces $H^p_A$,  $0<p\le 1$, defined by $L^2$-atoms
can be identified as the Triebel-Lizorkin spaces $F^0_{p2}$,
i.e. the Hardy spaces can be characterized via Littlewood-Paley square functions.
The same characterization of Hardy spaces in the setting of this article
can be obtain by using the method from \cite{DKKP}.
We will not pursue this line here.

Characterization of atomic Hardy spaces via other square functions
as well as their molecular decompositions are obtained
in \cite{HLMMY} (for $H^1$) and in \cite{Duong-Li} (for $H^p$, $0<p\le 1$)
in somewhat different settings. These can easily be adapted to our setting.
We will not elaborate on these results here.

The duality of atomic Hardy spaces and appropriately defined BMO and Lipschitz spaces
is established in \cite{HLMMY} (for $H^1$) and in \cite{Duong-Li} (for $H^p$, $0<p\le 1$)
in the settings of these articles.
The adaptation of these results to our setting is possible but is beyond the aims of this paper.

\section{Appendix}\label{sec:appendix-1}
\setcounter{equation}{0}

\subsection{Proof of Proposition~\ref{prop:dist}}

%
%
For the given $f\in \cS'$ there exist constants $m\in \bZ_+$ and $c>0$ such that (\ref{distribution-1}) holds.
Let $\phi\in\cS$. We have
$$
\varphi(\sqrt L)\phi(x) = \int_M \varphi(\sqrt L)(x,y)\phi(y) d\mu(y), \quad x\in M.
$$
To prove (\ref{phif}) we will interpret the above integral as a Bochner integral
over the Banach space 
$
V_m:=\{g\in \cap_{0\le \nu \le m} D(L^\nu): \|g\|_{V_m}:=\cP_m(g)<\infty\}
$
with $\cP_m$ defined in (\ref{norm-S}), see e.g. \cite{Yosida}, pp. 131-–133.
The completeness of $V_m$ follows (just as in the proof of \cite[Proposition 5.3]{KP})
by the fact that $L$ being a self-adjoint operator is also closed.
By the Hahn-Banach theorem the continuous linear functional $f$ can be extended to $V_m$ with the same norm.

Denote $F(y):=\varphi(\sqrt L)(\cdot,y)\phi(y)$. We have
\begin{align*}
\|F(y)\|_{V_m}
= \max_{0\le \nu\le m} \sup_{x\in M} (1+\rho(x, x_0))^m |L^\nu\varphi(\sqrt L)(x,y)\phi(y)|.
\end{align*}
By Theorem~\ref{thm:S-local-kernels},
applied with $f(\lambda)= \lambda^{2\nu}\varphi(\lambda)$,
it follows that $L^\nu\varphi(\sqrt{L})$
is an integral operator with a kernel satisfying the following inequality for any $\sigma >0$
$$
|L^\nu\varphi(\sqrt{L})(x, y)| \le c_\sigma |B(y, 1)|^{-1}(1+ \rho(x, y))^{-\sigma},
\quad 0\le \nu\le m.
$$
We choose $\sigma = m$.
On the other hand, as $\phi\in\cS$ we have, takeing into account (\ref{norm-S}),
$
|\phi(y)| \le \cP_\ell(\phi)(1+\rho(y, x_0))^{-\ell}
$
for any $\ell \ge 0$. We choose $\ell \ge m+2d+1$.
Putting these estimates together we get
$$
\|F(y)\|_{V_m}
\le c\max_{0\le \nu\le m} \sup_{x\in M}
\frac{\cP_\ell(\phi)(1+\rho(x, x_0))^m}{|B(y, 1)|(1+ \rho(x, y))^{m}(1+\rho(y, x_0))^{m+2d+1}}
$$
and using the obvious inequality $1+\rho(x, x_0) \le (1+ \rho(x, y))(1+\rho(y, x_0))$ we obtain
\begin{align*}
\|F(y)\|_{V_m}
&\le c \cP_\ell(\phi)|B(y, 1)|^{-1}(1+\rho(y, x_0))^{-2d-1}\\
&\le c \cP_\ell(\phi)|B(x_0, 1)|^{-1}(1+\rho(y, x_0))^{-d-1},
\end{align*}
where for the last inequality we used (\ref{D2}).
From the above and (\ref{tech-1}) it follows that
$
\int_M \|F(y)\|_{V_m} d\mu(y) \le c\cP_\ell(\phi).
$
Now, applying the theory of Bochner's integral we infer
$$
\Big\langle f, \int_M \varphi(\sqrt L)(\cdot,y)\phi(y)d\mu(y) \Big\rangle
= \int_M  \big\langle f, \varphi(\sqrt L)(\cdot,y) \big\rangle \overline{\phi(y)} d\mu(y).
$$
This coupled with (\ref{phi-distr}) implies (\ref{phif}).


We next prove (\ref{lipsch}); the proof of (\ref{slowly}) is simpler and will be omitted.
By the fact that (\ref{distribution-1}) holds for the given $f$ for some constants $m\in\bZ_+$ and $c>0$
and using (\ref{phif}) we obtain, for $x, x'\in M$,
\begin{equation}\label{phif-phif}
\begin{aligned}
&|\varphi(\sqrt L)f(x) - \varphi(\sqrt L)f(x')|
= |\langle  f , \varphi(\sqrt L)(x,\cdot)-\varphi(\sqrt L)(x',\cdot)\rangle\\
& \qquad\qquad
\le c\cP_m\big(\varphi(\sqrt L)(x,\cdot)-\varphi(\sqrt L)(x',\cdot)\big)\\
& \qquad\qquad
\le \max_{0\le \nu\le m} \sup_{y\in M}
(1+\rho(y, x_0))^m |L^\nu\varphi(\sqrt L)(x,y)-L^\nu\varphi(\sqrt L)(x',y)|.
\end{aligned}
\end{equation}
As above by Theorem~\ref{thm:S-local-kernels},
applied with $f(\lambda)= \lambda^{2\nu}\varphi(\lambda)$, it follows
that for any $\sigma >0$ and $0\le \nu\le m$
$$
|L^\nu\varphi(\sqrt{L})(x, y)-L^\nu\varphi(\sqrt{L})(x', y)|
\le c_\sigma |B(x, 1)|^{-1}\rho(x, x')^\alpha (1+ \rho(x, y))^{-\sigma}
$$
provided $\rho(x, x') \le 1$.
We choose $\sigma=m$.
We insert the above in (\ref{phif-phif}) and arrive at (\ref{lipsch}).
$\qed$

\subsection{Proof of Proposition~\ref{prop:convergence}}
This proof 
relies on the following


\begin{lemma}\label{lem:S}
Let $\sigma >0$ and $N\ge\sigma+d+\alpha/2$ with $\alpha>0$ from $(\ref{lip})$.
Then there exists a constant $c>0$ such that
for any $\phi\in \cS$ and $x, y\in M$
\begin{equation}\label{holder-S}
|\phi(x)-\phi(y)|
\le c \rho(x, y)^\alpha \cP_N(\phi)
\big[(1+\rho(x, x_0))^{-\sigma} + (1+\rho(y, x_0))^{-\sigma}\big].
\end{equation}
\end{lemma}

\noindent
{\bf Proof.}
Choose $\ph_0\in C^{\infty}(\bR_+)$ so that $0\le \ph_0 \le 1$,
$\ph_0(\lambda)=1$ for $\lambda \in [0, 1]$, and $\supp \ph_0\subset [0, 2]$.
Let $\ph(\lambda):=\ph_0(\lambda)-\ph_0(2\lambda)$
and set $\ph_j(\lambda):=\ph(2^{-j}\lambda)$, $j\ge 1$.
Clearly,
$\sum_{j \ge 0}\ph_j(\lambda)=1$ for $\lambda \in \bR_+$
and hence
$\phi=\sum_{j=0}^\infty \ph_j(\sqrt{L})\phi$
for $\phi\in\cS$ with the convergence in $L^\infty$ (see \cite[Proposition 5.5]{KP}).
Therefore,
$$
\phi(x)-\phi(y)=\sum_{j=0}^\infty \big(\ph_j(\sqrt{L})\phi(x)-\ph_j(\sqrt{L})\phi(y)\big),
\quad \forall x, y\in M, \; \forall \phi\in\cS.
$$
For $j\ge 1$ we have
\begin{equation}\label{phi-phi}
\begin{aligned}
&
\ph_j(\sqrt{L})\phi(x)-\ph_j(\sqrt{L})\phi(y)
= L^{-N}\ph_j(\sqrt{L})L^N\phi(x)-L^{-N}\ph_j(\sqrt{L})L^N\phi(y)\\
&\qquad\qquad
= \int_M \big[L^{-N}\ph(2^{-j}\sqrt{L})(x, z) - L^{-N}\ph(2^{-j}\sqrt{L})(y,z)\big]
L^N\phi(z)d\mu(z).
\end{aligned}
\end{equation}
Let $\omega(\lambda):= \lambda^{-2N}\ph(\lambda)$.
Then
$L^{-N}\ph(2^{-j}\sqrt{L})=2^{-2jN}\omega(2^{-j}\sqrt{L})$.
Clearly, $\omega\in C^\infty$ and $\supp \omega \subset [2^{-1}, 2]$.
Hence by Theorem~\ref{thm:S-local-kernels} it follows that
there exists a~constant $c_\sigma>0$ such that
\begin{equation}\label{local-LN-ph-1}
\big|L^{-N}\ph(2^{-j}\sqrt{L})(x, z)\big|
\le \frac{c_\sigma 2^{-2jN}}
{|B(x, 2^{-j})|\big(1+2^j\rho(x, z)\big)^{\sigma+d}}
\quad \hbox{and}
\end{equation}
\begin{equation}\label{local-LN-ph-2}
\big|L^{-N}\ph(2^{-j}\sqrt{L})(x, z) - L^{-N}\ph(2^{-j}\sqrt{L})(y, z)\big|
\le \frac{c_\sigma 2^{-2jN}\big(2^j\rho(x,y)\big)^\alpha}
{|B(x, 2^{-j})|\big(1+2^j\rho(x, z)\big)^{\sigma+d}},
\end{equation}
whenever $\rho(x, y) \le 2^{-j}$.

Fix $\phi\in\cS$. Then by (\ref{norm-S})
$|L^N\phi(z)| \le \cP_N(\phi)(1+\rho(z, x_0))^{-N}$, $z\in M$.

Let $\rho(x, y) \le 2^{-j}$.
The above, (\ref{phi-phi}), and (\ref{local-LN-ph-2}) yield
\begin{align*}
&|\ph_j(\sqrt{L})\phi(x)-\ph_j(\sqrt{L})\phi(y)|\\
&\quad\le c 2^{-j(2N-\alpha)}\rho(x,y)^\alpha \cP_N(\phi)
\int_M \frac{d\mu(z)}{|B(x, 2^{-j})|(1+2^j\rho(x, z))^{\sigma+d}(1+\rho(z, x_0))^N}\\
&\quad\le c 2^{-j(2N-d-\alpha)}\rho(x,y)^\alpha \cP_N(\phi)
\int_M \frac{d\mu(z)}{|B(x, 1)|(1+\rho(x, z))^{\sigma+d}(1+\rho(z, x_0))^{\sigma+d}}\\
&\quad\le \frac{c 2^{-j(2N-d-\alpha)}\rho(x,y)^\alpha \cP_N(\phi)}
{\big(1+\rho(x, x_0)\big)^{\sigma}}.
\end{align*}
Here we used that $|B(x,1)| \le c_0 2^{jd}|B(x, 2^{-j})|$, see (\ref{doubling}),
$N\ge \sigma+d$, and 
(\ref{tech-2}).

Let $\rho(x, y) > 2^{-j}$.
Using (\ref{local-LN-ph-1}) and some of the ingredients from above we get
\begin{align*}
&\Big|\int_M L^{-N}\ph(2^{-j}\sqrt{L})(x, z)L^N\phi(z)d\mu(z)\Big|\\
&\quad\le \int_M \frac{c 2^{-2jN} \cP_N(\phi)d\mu(z)}
{|B(x, 2^{-j})|\big(1+2^j\rho(x, z)\big)^{\sigma+d}(1+\rho(z, x_0))^N}\\
&\quad\le c 2^{-j(2N-d-\alpha)}\rho(x,y)^\alpha \cP_N(\phi)
\int_M \frac{d\mu(z)}{|B(x, 1)|(1+\rho(x, z))^{\sigma+d}(1+\rho(z, x_0))^{\sigma+d}}\\
&\quad\le \frac{c 2^{-j(2N-d-\alpha)}\rho(x,y)^\alpha \cP_N(\phi)}
{\big(1+\rho(x, x_0)\big)^{\sigma}}.
\end{align*}
Similarly
\begin{align*}
\Big|\int_M L^{-N}\ph(2^{-j}\sqrt{L})(y, z)L^N\phi(z)d\mu(z)\Big|
\le \frac{c 2^{-j(2N-d-\alpha)}\rho(x,y)^\alpha \cP_N(\phi)}
{\big(1+\rho(y, x_0)\big)^{\sigma}}.
\end{align*}
Putting the above estimates together we get for all $x, y\in M$ and $j\ge 1$
\begin{align}\label{est-ph-ph}
&|\ph_j(\sqrt{L})\phi(x)-\ph_j(\sqrt{L})\phi(y)|\\
&\qquad\le c 2^{-j(2N-d-\alpha)}\rho(x,y)^\alpha \cP_N(\phi)
\big[(1+\rho(x, x_0)\big)^{-\sigma} + (1+\rho(y, x_0)\big)^{-\sigma}\big]. \notag
\end{align}
In the same way, we use that (\ref{local-LN-ph-1})-(\ref{local-LN-ph-2})
hold for $\varphi_0(\sqrt{L})$ with $N=0$ to obtain
\begin{align*}
|\varphi_0(\sqrt{L})\phi(x)-\varphi_0(\sqrt{L})\phi(y)|
\le c\rho(x,y)^\alpha \cP_N(\phi)
\big[(1+\rho(x, x_0)\big)^{-\sigma} + (1+\rho(y, x_0)\big)^{-\sigma}\big].
\end{align*}
Summing up this estimate along with the estimates from (\ref{est-ph-ph})
($2N >d+\alpha$) we arrive at (\ref{holder-S}).
$\qed$

\smallskip

We are now in a position to prove Proposition~\ref{prop:convergence}.
Let  $\varphi \in \cS(\R)$, $\varphi$ be real-valued and even, and $\varphi(0)=1$.
It suffices to prove (\ref{decomp-dist-1}) only. Then (\ref{decomp-dist-2}) follows by duality,
see (\ref{phi-distr}).

Let $m\ge 0$ and $\phi\in\cS$.
Choose $\sigma >m+d+\alpha$ and $N\ge\sigma+d+\alpha/2$, where $\alpha>0$ is from $(\ref{lip})$.
By Theorem~\ref{thm:S-local-kernels}
$
|\varphi(\delta\sqrt{L})(x, y)|
\le c_\sigma |B(x, \delta)|^{-1}(1+\delta^{-1}\rho(x, y))^{-\sigma}
$
and
$\int_M \varphi(\delta\sqrt{L})(x, y)d\mu(y) = \varphi(0)=1$.
Therefore,
\begin{align*}
&(1+\rho(x, x_0))^m |L^m[\phi - \varphi(\delta\sqrt{L})\phi](x)|\\
&\qquad\qquad\qquad
= (1+\rho(x, x_0))^m \Big|\int_M \varphi(\delta\sqrt{L})(x, y)[L^m\phi(x) - L^m\phi(y)]d\mu(y)\Big|\\
&\qquad\qquad\qquad
\le c_\sigma(1+\rho(x, x_0))^m \int_M
\frac{|L^m\phi(x) - L^m\phi(y)|}{|B(x, \delta)|(1+\delta^{-1}\rho(x, y))^\sigma} d\mu(y)\\
&\qquad\qquad\qquad
= c_\sigma(1+\rho(x, x_0))^m \Big(\int_{B(x, 1)} \cdots +\int_{M\setminus B(x, 1)} \cdots \Big).
\end{align*}
As $\phi\in\cS$, then $L^m\phi\in\cS$ and applying Lemma~\ref{lem:S}
we obtain
\begin{align*}
&(1+\rho(x, x_0))^m \int_{B(x, 1)}
\frac{|L^m\phi(x) - L^m\phi(y)|}{|B(x, \delta)|(1+\delta^{-1}\rho(x, y))^\sigma} d\mu(y)\\
&\qquad\qquad
\le c(1+\rho(x, x_0))^m \int_{B(x, 1)}
\frac{\rho(x, y)^\alpha \cP_{m+N}(\phi)}
{|B(x, \delta)|(1+\delta^{-1}\rho(x, y))^\sigma(1+\rho(x, x_0))^\sigma} d\mu(y)\\
& \qquad\qquad
+ c(1+\rho(x, x_0))^m \int_{B(x, 1)}
\frac{\rho(x, y)^\alpha \cP_{m+N}(\phi)}
{|B(x, \delta)|(1+\delta^{-1}\rho(x, y))^\sigma(1+\rho(y, x_0))^\sigma} d\mu(y)\\
& \qquad\qquad
=: I_1+I_2.
\end{align*}
We use that $\sigma \ge m$, $\sigma-\alpha>d$, and (\ref{tech-1}) to obtain
\begin{align*}
I_1 &\le  c\cP_{m+N}(\phi)\int_{B(x, 1)}
\frac{\rho(x, y)^\alpha }{|B(x, \delta)|(1+\delta^{-1}\rho(x, y))^\sigma} d\mu(y)\\
& \le  c\cP_{m+N}(\phi)\int_{M}
\frac{\delta^\alpha }{|B(x, \delta)|(1+\delta^{-1}\rho(x, y))^{\sigma-\alpha}} d\mu(y)
\le  c\delta^\alpha\cP_{m+N}(\phi).
\end{align*}
Evidently, $1+\rho(x, x_0) \le (1+\rho(x, y))(1+\rho(y, x_0))$ and assuming $\delta\le 1$ we obtain
\begin{align*}
I_2 &\le  c\cP_{m+N}(\phi)\int_{B(x, 1)}
\frac{\rho(x, y)^\alpha }{|B(x, \delta)|(1+\delta^{-1}\rho(x, y))^{\sigma-m}} d\mu(y)\\
&\le  c\cP_{m+N}(\phi)\int_{M}
\frac{\delta^\alpha }{|B(x, \delta)|(1+\delta^{-1}\rho(x, y))^{\sigma-m-\alpha}} d\mu(y)
\le  c\delta^\alpha\cP_{m+N}(\phi).
\end{align*}
Here we also used that $\sigma>m+d+\alpha$ and (\ref{tech-1}).
Therefore, for any $x\in M$
\begin{equation}\label{est-B1}
(1+\rho(x, x_0))^m \int_{B(x, 1)}
\frac{|L^m\phi(x) - L^m\phi(y)|}{|B(x, \delta)|(1+\delta^{-1}\rho(x, y))^\sigma} d\mu(y)
\le c\delta^\alpha\cP_{m+N}(\phi).
\end{equation}

Since $\phi\in\cS$ we have by (\ref{norm-S})
$|L^m\phi(z)| \le \cP_{m+N}(\phi)(1+\rho(z, x_0))^{-N}$, $\forall z\in M$.
This leads to
\begin{align*}
&(1+\rho(x, x_0))^m \int_{M\setminus B(x, 1)}
\frac{|L^m\phi(x) - L^m\phi(y)|}{|B(x, \delta)|(1+\delta^{-1}\rho(x, y))^\sigma} d\mu(y)\\
&\le c\cP_{m+N}(\phi)\int_{M\setminus B(x, 1)}
\frac{(1+\rho(x, x_0))^m}{|B(x, \delta)|(1+\delta^{-1}\rho(x, y))^\sigma(1+\rho(x, x_0))^N} d\mu(y)\\
&+ c\cP_{m+N}(\phi)\int_{M\setminus B(x, 1)}
\frac{(1+\rho(x, x_0))^m}{|B(x, \delta)|(1+\delta^{-1}\rho(x, y))^\sigma(1+\rho(y, x_0))^N} d\mu(y)\\
&= J_1+J_2.
\end{align*}
Using that $N >\sigma > m$, $\sigma>d+\alpha$, (\ref{tech-1}),
and $\rho(x, y) \ge 1$ for $y\in M\setminus B(x, 1)$,
we get
\begin{align*}
J_1
&\le c\cP_{m+N}(\phi)\int_{M\setminus B(x, 1)}
\frac{d\mu(y)}{|B(x, \delta)|(1+\delta^{-1}\rho(x, y))^{\sigma}}\\
&\le c\cP_{m+N}(\phi)\int_{M}
\frac{\delta^\alpha d\mu(y)}{|B(x, \delta)|(1+\delta^{-1}\rho(x, y))^{\sigma-\alpha}}
\le c\delta^\alpha\cP_{m+N}(\phi).
\end{align*}
To estimate $J_2$ we use again that
$1+\rho(x, x_0) \le (1+\rho(x, y))(1+\rho(y, x_0))$ and assuming $\delta\le 1$ we obtain
\begin{align*}
J_2
&\le c\cP_{m+N}(\phi)\int_{M\setminus B(x, 1)}
\frac{d\mu(y)}{|B(x, \delta)|(1+\delta^{-1}\rho(x, y))^{\sigma-m}}\\
&\le c\cP_{m+N}(\phi)\int_{M}
\frac{\delta^\alpha d\mu(y)}{|B(x, \delta)|(1+\delta^{-1}\rho(x, y))^{\sigma-m-\alpha}}
\le c\delta^\alpha\cP_{m+N}(\phi).
\end{align*}
Consequently,
$$ 
(1+\rho(x, x_0))^m \int_{M\setminus B(x, 1)}
\frac{|L^m\phi(x) - L^m\phi(y)|}{|B(x, \delta)|(1+\delta^{-1}\rho(x, y))^\sigma} d\mu(y)
\le c\delta^\alpha\cP_{m+N}(\phi).
$$ 
This coupled with (\ref{est-B1}) leads to
$$
\sup_{x\in M}(1+\rho(x, x_0))^m |L^m[\phi - \varphi(\delta\sqrt{L})\phi](x)|
\le c\delta^\alpha\cP_{m+N}(\phi),
$$
which yields (\ref{decomp-dist-1}).

The proof of (\ref{decomp-dist-2}) in $L^p$ for $f\in L^p$ is straightforward
and will be omitted.

The almost everywhere convergence
$\lim_{\ddelta\to 0} \varphi (\ddelta\sqrt L) f(x)=f(x)$
for $f\in L^p(M)$, $1\le p\le \infty$,
follows by a standard argument using the doubling condition (\ref{doubling-0}),
the weak $(1, 1)$ boundedness of the Hardy-Littlewood maximal operator,
and the nearly exponential localization of the summability kernel $\varphi (\ddelta\sqrt L)(x, y)$.
$\qed$

\subsection{Proof of Lemma \ref{lem:construction}.}

%
%
Suppose $\phi\in C^\infty_0(\R)$, $\phi\ge 0$, $\supp \phi\subset [-1/4, 1/4]$,
$\phi(\xi)>0$ for $\xi\in (-1/4, 1/4)$, and $\phi$ is even.
Let $\Theta(\xi):= \phi(\xi+1/2) - \phi(\xi-1/2)$ for $\xi\in \R$.
Clearly $\Theta$ is odd.

There is no loss of generality in assuming that $m$ is even,
for otherwise we work with $m+1$ instead.
Denote $\Delta_h^m:= (T_h-T_{-h})^m$, where
$T_hf(\xi):= f(\xi+h)$.
Define
$$
\varphi(x):= \int_\R \xi^{-1}\Delta_h^m\Theta (\xi) e^{i\xi x}d\xi
= 2\pi \cF^{-1}\big(\xi^{-1}\Delta_h^m\Theta (\xi)\big),
\quad x\in\R, \;\; h:=\frac{1}{8m}.
$$
Evidently, $\varphi\in \cS(\R)$, $\varphi$ is even and real-valued,
and $\hat\varphi(\xi) = 2\pi\xi^{-1}\Delta_h^m\Theta (\xi)$.
Hence $\supp \hat\varphi \subset [-1, 1]$.
Furthermore, for $\nu =1, 2, \dots, m$,
$$
\varphi^{(\nu)}(0)= \int_\R \xi^{\nu-1}\Delta_h^m\Theta (\xi)d\xi
= (-1)^m \int_\R \Theta (\xi)\Delta_h^m \xi^{\nu-1}d\xi =0
\quad\hbox{and}
$$
$$
\varphi(0)= \int_\R \xi^{-1}\Delta_h^m\Theta (\xi)d\xi
= (-1)^m \int_\R \Theta (\xi)\Delta_h^m \xi^{-1}d\xi
= 2(-1)^m \int_{1/4}^{3/4} \Theta (\xi)\Delta_h^m \xi^{-1}d\xi.
$$
However, for any sufficiently smooth function $f$ we have
$\Delta_h^mf(\xi) = (2h)^mf^{(m)}(\theta)$, where $\theta \in (\xi-mh, \xi+mh)$.
Hence,
$$
\Delta_h^m \xi^{-1} = (2h)^m m!(-1)^m \theta^{-m-1}
\quad\hbox{with}\quad \theta \in (\xi-mh, \xi+mh)\subset [1/8, 7/8].
$$
Therefore, $\varphi(0)\ne 0$ and then $\varphi(0)^{-1}\varphi(x)$
has the claimed properties.
$\qed$

\subsection{Proof of Lemma \ref{lem:Whitney}.}
%
%
Choose $\{B(\xi_j, \rho(\xi_j)/5)\}_{j\in\bN}$
to be a maximal disjoint subcollection of $\{B\big(x, \rho(x)/5\big)\}_{x\in\Omega}$,
whose existence follows by Zorn's lemma.
Then (b) is obvious.


We now establish (a). Assume to the contrary that there exists $x\in \Omega$
such that $x\not\in\cup_{j\in \bN} B(\xi_j, \rho_j/2)$.
From the construction of $\{B(\xi, \rho_j/5)\}_{j\in\bN}$ it follows that
$B(x, \rho(x)/5)\cap B(\xi, \rho_j/5)\ne \emptyset$ for some $j\in\bN$.
We claim that
\begin{equation}\label{est-rho}
\rho(\xi_j) > (2/3)\rho(x).
\end{equation}
Indeed, assume that $\rho(\xi_j) \le (2/3)\rho(x)$.
Then
$$
\rho(x, \xi_j) < (1/5)(\rho(\xi_j)+\rho(x)) \le (1/3)\rho(x).
$$
Therefore,
$
B(\xi_j, \rho_j) \subset B\big(x, \rho(x, \xi_j)+ \rho(\xi_j)\big)
\subset B(x, \rho(x)),
$
where the first inclusion is strict.
This implies
$B\big(\xi_j, (1+\eta)\rho_j\big) \subset B(x, \rho(x)) \subset \Omega$
for some $\eta>0$.
But from the definition of $\rho_j$ it follows that
$B\big(\xi_j, (1+\eta)\rho_j\big)\cap \Omega^c \ne \emptyset$.
This is a contradiction which proves (\ref{est-rho}).
From (\ref{est-rho}) we infer
$$
\rho(x, \xi_j) < (1/5)(\rho(\xi_j)+\rho(x)) \le (1/5)(1+3/2)\rho(\xi_j)
=(1/2)\rho(\xi_j),
$$
which verifies (a).


To prove (c) assume
$B\big(\xi_j, \frac{3\rho_j}{4}\big)\cap B\big(\xi_\nu, \frac{3\rho_\nu}{4}\big)\ne \emptyset$
for some $j, \nu\in\bN$.
We will show that $\rho_j \le 7\rho_\nu$.
We proceed similarly as above.
Assume that $\rho_j > 7\rho_\nu$.
Then $\rho(\xi_j, \xi_\nu) \le (3/4)(\rho_j+\rho_\nu) \le (6/7)\rho_j$
yielding
$$
B(\xi_\nu, \rho_\nu) \subset B\big(\xi_j, \rho(\xi_j, \xi_\nu)+ \rho_\nu\big)
\subset B\big(\xi_j, (6/7)\rho_j+ (1/7)\rho_j\big) = B(\xi_j, \rho_j),
$$
where the first inclusion is strict.
As above this leads to a contradiction which shows that  $\rho_j \le 7\rho_\nu$.


To prove (d),
assume that balls $B(\xi_{\nu_m}, 3\rho_{\nu_m}/4)$, $m=1, 2, \dots, K$, intersect $B(\xi_j, 3\rho_j/4)$.
Then from above $\rho_j \le 7\rho_{\nu_m}$, $m=1, 2, \dots, K$.
Using this, (\ref{D2}) and (\ref{doubling}) we get
\begin{align*}
|B(\xi_j, 8\rho_j)| &\le c_0\Big(1+\frac{\rho(\xi_j, \xi_{\nu_m})}{8\rho_j}\Big)^d|B(\xi_{\nu_m}, 8\rho_j)|\\
&\le c_0^2\Big(1+\frac{\rho(\xi_j, \xi_{\nu_m})}{8\rho_j}\Big)^d
40^d |B(\xi_j, \rho_{\nu_m}/5)|.
\end{align*}
However, using (c),
$\rho(\xi_j, \xi_{\nu_m}) \le (3/4)(\rho_j+\rho_{\nu_m}) \le 6\rho_j$.
Therefore,
$$
|B(\xi_j, 8\rho_j)| \le c_0^2 70^d |B(\xi_j, \rho_{\nu_m}/5)|
$$
and summing up we obtain
\begin{equation}\label{B-B}
K|B(\xi_j, 8\rho_j)| \le 70^d c_0^2 \sum_{m=1}^K |B(\xi_j, \rho_{\nu_m}/5)|.
\end{equation}
On the other hand, by (b) the balls $B(\xi_{\nu_m}, \rho_{\nu_m}/5)$, $m=1, \dots, K$, are disjoint,
and since each ball $B(\xi_{\nu_m}, 3\rho_{\nu_m}/4)$ intersects $B(\xi_j, 3\rho_j/4)$
and $\rho_{\nu_m} \le 7\rho_j$ we have
$$
B(\xi_{\nu_m}, \rho_{\nu_m}/5) \subset B\big(\xi_j, 3\rho_j/4+ (3/4+1/5)\rho_{\nu_m}\big)
\subset B(\xi_j, 8\rho_j).
$$
Consequently,
$
\sum_{m=1}^K |B(\xi_{\nu_m}, \rho_{\nu_m}/5)| \le |B(\xi_j, 8\rho_j)|.
$
This coupled with (\ref{B-B}) yields
$K \le 70^d c_0^2$.
$\qed$

\end{document}